\documentclass{amsart}%
\usepackage{amsthm}
\usepackage{amsmath}
\usepackage{amsfonts}
\usepackage{amssymb}
\usepackage{dsfont}
\usepackage{esint}
\usepackage{graphicx}
\usepackage{color}
\usepackage{url}
\usepackage[latin1]{inputenc}
\usepackage[T1]{fontenc}
\usepackage[english]{babel}%
\numberwithin{equation}{section}

\usepackage{amsaddr}

\providecommand{\U}[1]{\protect\rule{.1in}{.1in}}

\newtheorem{theorem}{Theorem}[section]

\newtheorem{definition}{Definition}[section]

\newtheorem{lemma}{Lemma}[section]

\newtheorem{proposition}{Proposition}[section]
\newtheorem{remark}{Remark}[section]

\newcommand{\rst}[1]{\ensuremath{\raise-1.0ex\hbox{\large{$\vert_{#1}$}}}}

\begin{document}

\title[Existence for a mixture of two power-law fluids]{Existence for the steady problem of a mixture of two power-law fluids}

\date{\textbf{April 26, 2012}}
\author[H.B. de Oliveira]{Hermenegildo Borges de Oliveira$^{\ast,\ast\ast}$}

\email{holivei@ualg.pt}

\address{$^{\ast}$FCT - Universidade do Algarve and $^{\ast\ast}$CMAF - Universidade de Lisboa, Portugal.}

\thanks{The authors work was partially supported by FEDER and FCT-Plu\-ri\-an\-ual 2010 and by the grant SFRH/BSAB/1058/2010, MCTES, Portugal.}

\begin{abstract}
The steady problem resulting from a mixture of two distinct fluids of power-law type is analyzed in this work. %
Mathematically, the problem results from the superposition of two power laws, one for a constant power-law index with other for a variable one. %
For the associated boundary-value problem, we prove the existence of very weak solutions, provided the variable power-law index is bounded from above by the constant one. %
This result requires the lowest possible assumptions on the variable power-law index and, as a particular case, extends the existence result by Ladyzhenskaya~\cite{L-1967}  to the case of a variable exponent and for all zones of the pseudoplastic region. %
In a distinct result, we extend a classical theorem on the existence of weak solutions to the case of our problem.



\end{abstract}

\maketitle

\maketitle

\selectlanguage{english}


%\author{}

%\date{}

%\date{}
%\dedicatory{To my parents}
%%% ----------------------------------------------------------------------

\section{Introduction}\label{Sect-Int}

In Fluid Mechanics the constitutive equation relates the stress on a fluid element to other fluid quantities by the relation
\begin{equation*}
\mathbf{T}=-p\mathbf{I}+\mathbf{S}\,,
\end{equation*}
where $\mathbf{T}$ is the Cauchy stress tensor, $p$ stands for the pressure, $\mathbf{I}$ is the unit tensor, $\mu$ is the dynamic viscosity and $\mathbf{S}$ is the deviatoric part of the stress tensor.
The simplest model of Fluid Mechanics is the Newtonian fluid which by definition is a fluid that continues to flow, regardless of the forces acting on it. %
For these fluids, the constitutive equation is the Stokes law and therefore the deviatoric part of the stress tensor is given by
\begin{equation*}
\mathbf{S}=2\mu\mathbf{D}\,,\quad \mathbf{D}\equiv\mathbf{D}(\mathbf{u})=\frac{1}{2}\left(\mathbf{\nabla\,u}+\mathbf{\nabla\,u}^T\right)\,,
\end{equation*}
where $\mathbf{D}$ is the rate of strain tensor and $\mathbf{u}$ is the velocity field. %
Although the Newtonian fluid model captures the characteristics of many fluids as water solutions, gasoline, vegetal and mineral oils, this model is quite inadequate for describing the complex rheological behavior of many other fluids. In this case are all the fluids in which the apparent fluid viscosity $\mu$ decreases or increases with the shear rate history $|\mathbf{D}|$. These fluids form a large class called non-Newtonian, or generalized Newtonian fluids, and are divided into pseudoplastic and dilatant fluids. In pseudoplastic fluids the viscosity gradually decreases with increasing shear rate and, due to this, they are often called shear thinning fluids. Shear thickening fluids is the other name found in the literature for dilatant fluids and are fluids in which the viscosity increases with the shear rate. Examples of pseudoplastic fluids are milk fluids, varnishes, shampoo and blood fluids, whereas polar ice, glaciers, volcano lava and sand are all examples of dilatant fluids. %
There are also some another class, called Bingham fluids, that are similar to pseudoplastic fluids, but they exhibit a yield point.
Examples of this fluids are drilling muds used in petroleum industry, toothpaste and face creams. %
The existence of a yield point means fluid flow is prevented below a critical stress level, but flow occurs when the critical stress level is exceeded.
During the 20th century a large number of models have been proposed in the literature to model all types of non-Newtonian
fluids under diverse flow conditions. However, it was only since the pioneer work done by Oldroyd~\cite{O-1950}, during the 1950's, that were established some guiding principles  to provide a constructive framework for the formulation of constitutive equations to the Cauchy stress tensor (see Barnes \emph{et al.}~\cite{BHW-1993}). %
The power-law, or Ostwald-de Waele model, is one of the simplest non-Newtonian fluid models then introduced and can be defined by the relation
\begin{equation}\label{p-l-const}
\mathbf{S}=\mu|\mathbf{D}|^{\gamma-2}\mathbf{D}\,,
\end{equation}
where $\gamma\geq 1$ characterizes the flow behavior and is usually called the power-law index. %
The power-law is often used to model pseudoplastic fluids though it can also be used for modeling dilatant fluids or for Newtonian fluids (see \emph{e.g.} Schowalter~\cite{S-1978}):
$$
\mbox{power-law}\
                  \left\{\begin{array}{ll}
                         \mbox{Bingham} & \mbox{if $\gamma=1$} \\
                         \mbox{pseudoplastic}\quad & \mbox{if $1<\gamma<2$} \\
                         \mbox{Newtonian} & \mbox{if $\gamma=2$} \\
                         \mbox{dilatant} & \mbox{if $\gamma>2$}\,.
                        \end{array}\right.
$$
There are also some fluids that cannot be cataloged into a single class of non-Newtonian fluids. These fluids can go, for instance, from the consistency of a liquid to that of a gel, and back, with response times on the order of milliseconds. In this case are the electrorheological fluids whose rheological properties are controllable through the application of an electric field, showing useful and special function with the effect of reversibility. The best example of electrorheological fluids are suspensions dispersed with some polymeric colloids which show trembling shear behavior under an electric field (see \emph{e.g.} Hao~\cite{H-2005}). Due to this, in the sequel, by a \emph{trembling fluid} we mean a fluid with a variable power-law index, whereas \emph{sustaining fluid} is the name we shall use for a fluid with a constant power-law index. In this way, trembling fluids can also be medelled by the power-law (\ref{p-l-const}), but with the significant difference that now the power-law index $\gamma$ may vary with other quantity under study or, in the simplest cases, with the space and time position (see \emph{e.g.} Rajagopal and R\.{u}\v{z}i\v{c}ka~\cite{RR-2001}). %

In this article we will study the mathematical problem for a steady motion of a generalized fluid contained in a bounded domain $\Omega\subset\mathds{R}^N$, $N\geq 2$, with the boundary denoted by $\partial\Omega$. We assume the motion is described by the following boundary-value problem for the generalized Navier-Stokes equations:
\begin{equation}\label{geq1-inc}
\mathrm{div}\,\mathbf{u}=0\quad\mbox{in}\quad \Omega;
\end{equation}
\begin{equation}\label{geq1-vel}
\mathbf{div}(\mathbf{u}\otimes\mathbf{u})
=\mathbf{f}-\mathbf{\nabla}p+\mathbf{div}\,\mathbf{S} \quad\mbox{in}\quad \Omega;
\end{equation}
\begin{equation}\label{geq1-bc}
\mathbf{u}=\mathbf{0} \qquad\mbox{on}\quad \partial\Omega;
\end{equation}
where, $p$ stands, now, for the pressure divided by
the constant density and $\mathbf{f}$ is the external forces field. %
We assume the dependence of the deviatoric stress tensor $\mathbf{S}$ on the space variable $\mathbf{x}$ and on the strain rate tensor $\mathbf{D}$ is given by the superposition of two different power laws, one for a \emph{sustaining} fluid with other for a \emph{trembling} one:
\begin{equation}\label{e-stress}
\mathbf{S}=\left(\mu_1|\mathbf{D}|^{\gamma-2}+\mu_2|\mathbf{D}|^{q(\mathbf{x})-2}\right)\mathbf{D}\,.
\end{equation}
Here, $\mu_1$ and $\mu_2$ are positive constants related with the fluid viscosity, and $\gamma$ and $q$ are the sustaining and trembling power-law indexes, respectively: $\gamma$ is a constant  and $q$ depends on the space variable. %
The specific physical problem we are interested in and where (\ref{e-stress}) can be potentially useful is the case of a dilute suspension of an electrorheological material in a dilatant fluid. %
The object of superposition of generalized fluids is to produce flow patterns similar to those of practical interest. %
Many systems, among them polymer solutions and emulsions, behave, in the dilute regime, as superposed fluids in the above sense.
The best example are polymer solutions in which the polymer segments tend to repel each other, since they prefer contact the solvent molecules rather then among themselves (see \emph{e.g.} Oswald~\cite{O-2009}).
Moreover, superposition of fluids is justified, in the light of theoretical mechanics, as a powerful tool to replace the Boltzman superposition principle in the case of materials with nonlinear behavior (see \emph{e.g.} Dealy~\cite{D-2009}). %
The constitutive relation (\ref{e-stress}) can be used to model many other generalized fluids as follows:
$$
\left\{
\begin{array}{ll}
    \mbox{sustaining power-law} & \mbox{if  $\mu_1>0$ and $\mu_2=0$} \\
    \mbox{trembling power-law}\quad & \mbox{if  $\mu_1=0$ and $\mu_2>0$\,.}
  \end{array}\right.
$$
In particular, by making $\gamma=2$ in (\ref{e-stress}), we obtain a generalization of the Sisko model to the trembling fluids.
The Sisko model  has  been  checked  experimentally  to  fit  accurately  the  viscosity  data  of  various  commercial  greases  made  from
petroleum  oils with  one  of the standard  thickening  agents  such  as calcium  fatty  acid, lithium
hydroxy  stearate,  sodium  tallow  or  hydrophobic  silica  over  a  wide  range  of  shear  rate (see Sisko~\cite{S-1958}). %
Letting also $\gamma=2$, we recover a specific Carreau-Yassuda model which is very often used to describe blood flows (see \emph{e.g.} Carreau \emph{et al.}~\cite{CD-HC-1997}):
\begin{equation*}\label{CY-model}
\mathbf{S}=\mu_{\infty}+(\mu_{0}-\mu_{\infty})\left(1+\lambda|\mathbf{D}|^{a}\right)|^{\frac{n-2}{a}}\,.
\end{equation*}
In the specific example we want to address, $\mu_0$ stands for the zero-shear viscosity, the infinite-shear viscosity $\mu_{\infty}$ is zero, where $\lambda$ is a relaxation time, the power-law index is absent, \emph{i.e.} $n=2$, and $a$ stands for a variable shape parameter. %

The outline of this work is the following. %
In Section~\ref{Sect-Int}, not only we presented the problem that we shall study in this work, but also we have given a physical motivation for doing so. %
The main notation used throughout the text and some auxiliary results are presented in Section~\ref{S-Prel}. %
Section~\ref{Sect-HB} is devoted to review the main existence results for some particular cases of the problem (\ref{geq1-inc})-(\ref{e-stress}).
In Section~\ref{Sect-WF} we define the notion of solutions we shall consider and we state two different existence results: Theorems~\ref{th0-exst-ws-pp} and~\ref{th-exst-ws-pp}. %
From Section~\ref{Sect-Exist-RP} to Section~\ref{Sect-Conv-EST} we shall prove the main result of this work: Theorem~\ref{th-exst-ws-pp}. %

\section{Preliminaries}\label{S-Prel}

The notation used in this work is largely standard in Mathematical Fluid Mechanics (see \emph{e.g.} Lions~\cite{Lions-1969}).
In this article, the notations $\Omega$ or $\omega$ stand always for a domain, \emph{i.e.}, a connected open subset of $\mathds{R}^N$, $N\geq 1$.
Given $k\in\mathds{N}$, we denote by $\mathrm{C}^{k}(\Omega)$ the space of all $k$-differentiable functions in $\Omega$. %
By $\mathrm{C}^{\infty}_0(\Omega)$ we denote the space of all infinity-differentiable functions with compact support in $\Omega$. %
In the context of distributions, the space $\mathrm{C}^{\infty}_0(\Omega)$ is denoted by $\mathcal{D}(\Omega)$ instead. %
The space of distributions over $\mathcal{D}(\Omega)$ is denoted by $\mathcal{D}'(\Omega)$. %
If $\mathrm{X}$ is a generic Banach space, its dual space is denoted by $\mathrm{X}'$.
Let $1\leq q\leq \infty$ and $\Omega\subset\mathds{R}^N$, with $N\geq 1$, be a domain. %
We use the classical Lebesgue spaces $\mathrm{L}^q(\Omega)$, whose norm is
denoted by $\|\cdot\|_{\mathrm{L}^q(\Omega)}$. %
For any nonnegative $k$,
$\mathrm{W}^{k,q}(\Omega)$ denotes the Sobolev space of all
functions $u\in\mathrm{L}^q(\Omega)$ such
that the weak derivatives $\mathrm{D}^{\alpha}u$ exist, in the generalized sense, and are in
$\mathrm{L}^q(\Omega)$ for any multi-index $\alpha$ such that
$0\leq |\alpha|\leq k$.
In particular, $\mathrm{W}^{1,\infty}(\Omega)$ stands for the space of Lipschitz functions. %
The norm in $\mathrm{W}^{k,q}(\Omega)$ is denoted by
$\|\cdot\|_{\mathrm{W}^{k,q}(\Omega)}$.
We define $\mathrm{W}^{k,q}_0(\Omega)$ as the closure of $\mathrm{C}^{\infty}_0(\Omega)$ in $\mathrm{W}^{k,q}(\Omega)$.
For the dual space of $\mathrm{W}^{k,q}_0(\Omega)$, we use the identity $(\mathrm{W}^{k,q}_0(\Omega))'=\mathrm{W}^{-k,q'}(\Omega)$, up to an isometric isomorphism. %
We shall distinguish (second-order) tensor-valued and vector-valued space functions from scalar-valued ones by using boldface letters.
Although we are going to use the same notation for tensor-valued and vector-valued space functions, the distinction between these spaces will always be clear from the exposition.

We denote by $\mathcal{P}(\Omega)$ the set of all measurable functions $q:\Omega\to [1,\infty]$ and define
$$\displaystyle q^{-}:=\mathrm{ess}\inf_{\hspace{-0.5cm}x\in \Omega}q(x),\quad q^{+}:=\mathrm{ess}\sup_{\hspace{-0.5cm}x\in \Omega}q(x).$$
Given $q\in\mathcal{P}(\Omega)$, we denote by $\mathrm{L}^{q(\cdot)}(\Omega)$ the space of all measurable functions $f$ in $\Omega$ such that its semimodular is finite:
\begin{equation}\label{ap3}
A_{q(\cdot)}(f):=\int_{\Omega}|f(x)\mathbf{|}^{q(x)}d\,x<\infty.
\end{equation}
The space $\mathrm{L}^{q(\cdot)}(\Omega)$ is called Orlicz space and is also known by Lebesgue space with variable exponent. %
Equipped with the norm
\begin{equation}\label{ap4}
\left\Vert f\right\Vert_{\mathrm{L}^{q(\cdot)}(\Omega)}:=
\inf\left\{\kappa>0 : A_{q(\cdot)}\left(\frac{f}{\kappa}\right) \leq
1\right\},
\end{equation}
$\mathrm{L}^{q(\cdot)}(\Omega)$ becomes a Banach space. %
Note that the infimum in (\ref{ap4}) is attained if $A_{q(\cdot)}(f)>0$. %
If $q^{+}<\infty$, $\mathrm{L}^{q(\cdot)}(\Omega)$ is separable and the space $\mathrm{C}^{\infty}_0(\Omega)$ is dense in $\mathrm{L}^{q(\cdot)}(\Omega)$.
Moreover, if
\begin{equation}\label{ap1-A}
1<q^{-}\leq q^{+}<\infty,
\end{equation}
$\mathrm{L}^{q(\cdot)}(\Omega)$ is reflexive. %
One problem in Orlicz spaces, is the relation between the semimodular (\ref{ap3}) and the norm (\ref{ap4}). %
If (\ref{ap1-A}) is satisfied, one can shows that
\begin{equation}\label{ap4-56}
\|f\|_{\mathrm{L}^{q(\cdot)}(\Omega)}^{q^{-}}-1\leq A_{q(\cdot)}(f)\leq
\|f\|_{\mathrm{L}^{q(\cdot)}(\Omega)}^{q^{+}}+1\,.
\end{equation}
In Orlicz spaces, there holds a version of H\"{o}lder's inequality, called generalized H\"{o}lder's inequality.

Given $q\in\mathcal{P}(\Omega)$, the Orlicz-Sobolev space $W^{1,q(\cdot)}(\Omega)$ is defined as:
$$W^{1,q(\cdot)}(\Omega):=\left\{f\in \mathrm{L}^{q(\cdot)}(\Omega): \mathrm{D}^{\alpha}f\in\mathrm{\mathrm{L}}^{q(\cdot)}(\Omega),\ 0\leq |\alpha|\leq 1\right\}\,.$$
In the literature, this space is also known by Sobolev space with variable exponent. %
In $W^{1,q(\cdot)}(\Omega)$ is defined a semimodular and the correspondent induced norm analogously as in (\ref{ap3})-(\ref{ap4}), which are equivalent, respectively, to
\begin{equation*}\label{ap3-W}
A_{1,q(\cdot)}(f):=A_{q(\cdot)}(f)+A_{q(\cdot)}(|\mathbf{\nabla}f|)
\end{equation*}
and
\begin{equation*}\label{ap4-W}
\|f\|_{W^{1,q(\cdot)}(\Omega)}:=\|f\|_{\mathrm{L}^{q(\cdot)}(\Omega)}+\|\mathbf{\nabla}f\|_{\mathrm{L}^{q(\cdot)}(\Omega)}\,.
\end{equation*}
For this norm, $W^{1,q(\cdot)}(\Omega)$ is a Banach space, which becomes separable and reflexive in the same conditions as $\mathrm{L}^{q(\cdot)}(\Omega)$. %
The Orlicz-Sobolev space with zero boundary values is defined by:
$$W^{1,q(\cdot)}_0(\Omega):=\overline{\left\{f\in W^{1,q(\cdot)}(\Omega): \mathrm{supp}\ f\subset\subset \Omega\right\}}^{\ \|\cdot\|_{W^{1,q(\cdot)}(\Omega)}}\,.$$
In contrast to the case of classical Sobolev spaces, the set $\mathrm{C}_0^{\infty}(\Omega)$ is not necessarily dense in $\mathrm{W}^{1,q(\cdot)}_0(\Omega)$ -- the closure of $\mathrm{C}_0^{\infty}(\Omega)$ in $\mathrm{W}^{1,q(\cdot)}(\Omega)$ is strictly contained in $\mathrm{W}^{1,q(\cdot)}_0(\Omega)$. %
The equality holds only if $q$ is globally log-Hölder continuous, \emph{i.e.}, if exist positive constants
$C_1$, $C_2$ and $q_{\infty}$ such that
\begin{equation}\label{glob-log-H}
\left|q(\mathbf{x})-q(\mathbf{y})\right|\leq \frac{C_1}{\ln(e+1/|\mathbf{x}-\mathbf{y}|)},\quad
\left|q(\mathbf{x})-q_{\infty}\right|\leq \frac{C_2}{\ln(e+|\mathbf{x}|)}\quad\forall\ \mathbf{x},\ \mathbf{y}\in \Omega.
\end{equation}
For a thorough analysis on Orlicz and Orlicz-Sobolev spaces, we address the reader to the monograph by Diening \emph{et al.}~\cite{DHHR}.

\section{Historical background}\label{Sect-HB}

To the best of our knowledge, problem (\ref{geq1-inc})-(\ref{e-stress}) is new, but there is an extensive literature on some particular cases of it.
Here we shall be concerned only with the existence results for the problem (\ref{geq1-inc})-(\ref{e-stress}) or its simplifications. %

In the case of $\mu_1\not=0$ and $\mu_2=0$, the main achievements on the existence results for this problem were done by Lions~\cite{Lions-1969}, Frehse \emph{et al.}~\cite{FMS-1997,FMS-2003} and R\.{u}\-\v{z}i\-\v{c}ka~\cite{R-1997}. %
But even earlier than the first of these authors, the mathematical analysis of the problem (\ref{geq1-inc})-(\ref{e-stress}) with $\mu_1\not=0$, but $q\equiv 2$  and $N=3$, was performed by Ladyzhenskaya in the work~\cite{L-1967}. It should be noted that Ladyzhenskaya's problem~\cite{L-1967} is, in fact, the problem considered by Sisko~\cite{S-1958}, more or less 10 years before, in rheological experiments. Surprisingly is that, in the literature, there is no evidence that Ladyzhenskaya knew Sisko's work. %
Moreover, while Sisko was working on experiments with various commercial grease flows, Ladyzhenskaya studied this problem motivated by the important issue of uniqueness for the classical Navier-Stokes problem.
Problem (\ref{geq1-inc})-(\ref{e-stress}) with $\mu_1\not=0$ and $\mu_2=0$ was consider by Lions~\cite{Lions-1969} using the same approach of Ladyzhenskaya~\cite{L-1967}. In both works \cite{L-1967} and \cite{Lions-1969}, and under the assumption that $\mathbf{f}\in \mathbf{V}_\gamma'$ and
\begin{equation}\label{N-LL}
\gamma\geq\frac{3N}{N+2},
\end{equation}
the authors have proved the existence of weak solutions in the class $\mathbf{V}_\gamma$ satisfying to the energy relation
\begin{equation}\label{er-LL}
\int_{\Omega}\left(\mathbf{S}(\mathbf{D}(\mathbf{u}))-\mathbf{u}\otimes\mathbf{u}\right):\mathbf{D}(\mathbf{\varphi})\,d\mathbf{x}
=\int_{\Omega}\mathbf{f}\cdot\mathbf{\varphi}\,d\,\mathbf{x}
\end{equation}
for all $\mathbf{\varphi}\in\mathbf{V}_\gamma$. Here the notation $\mathbf{V}_\gamma'$ stands for the dual space of $\mathbf{V}_\gamma$ and
\begin{equation}\label{V-reg}
\mathcal{V}:=\{\mathbf{v}\in\mathbf{C}^{\infty}_0(\Omega):\mathrm{div\,}\mathbf{v}=0\},
\end{equation}
\begin{equation}\label{V-q=const}
\mathbf{V}_\gamma:=\mbox{closure of $\mathcal{V}$ in $\mathbf{W}^{1,\gamma}(\Omega)$}.
\end{equation}
The proofs in \cite{L-1967,Lions-1969} use the theory of monotone operators together with compactness arguments, and the lower bound (\ref{N-LL}) results from controlling the boundedness of the convective term  $\mathbf{u}\otimes\mathbf{u}:\mathbf{D}(\mathbf{\varphi})$ in $\mathbf{L}^1(\Omega)$ for $\mathbf{u}$ and $\mathbf{\varphi}$ belonging to $\mathbf{V}_\gamma$.
More or less 30 years later the results of \cite{L-1967,Lions-1969} were improved in~\cite{FMS-1997,R-1997} for values of $\gamma$ such that
\begin{equation}\label{N-FMS-R}
\gamma>\frac{2N}{N+1}.
\end{equation}
Assuming that $\mathbf{f}\in \mathbf{L}^{\gamma'}(\Omega)$ in \cite{FMS-1997} and $\mathbf{f}\in \mathbf{V}_\gamma'$ in \cite{R-1997},
the authors have established existence results in the same class $\mathbf{V}_\gamma$ of \cite{L-1967,Lions-1969} satisfying to the energy relation (\ref{er-LL}), but for all $\mathbf{\varphi}\in\mathcal{V}$. %
In order to distinguish the weak solutions found in \cite{L-1967,Lions-1969} from these established in \cite{FMS-1997,R-1997}, we shall denote the later as \emph{very weak solutions}. %
An important feature of the works \cite{FMS-1997,R-1997}, is that there the space of test functions allows us to consider test functions with more regularity. %
Due to that, to control the boundedness of $\mathbf{div}(\mathbf{u}\otimes\mathbf{u})\cdot\mathbf{\varphi}$ in $\mathbf{L}^1(\Omega)$ for $\mathbf{u}$ in $\mathbf{V}_\gamma$ and $\mathbf{\varphi}$ in $\mathcal{V}$, led us to a lower bound for $\gamma$ (\emph{f.} (\ref{N-FMS-R})). %
Although this important difference, the proofs of \cite{FMS-1997} and \cite{R-1997} also use the theory of monotone operators together with compactness arguments and they differ only on a slight different application of the $L^{\infty}$-truncation method. %
A little bit later, in~\cite{FMS-2003} the authors improved theirs own result~\cite{FMS-1997} for
\begin{equation}\label{N-FMS}
\gamma>\frac{2N}{N+2}.
\end{equation}
The new feature was the application of the Lipschitz-truncation method. %
By this method, the authors could use all the regularity they needed for the test function and to worry only about the boundedness control of $\mathbf{u}\otimes\mathbf{u}$ in $\mathbf{L}^1(\Omega)$ for $\mathbf{u}$ in $\mathbf{V}_\gamma$ and $\mathbf{\varphi}$ in $\mathcal{V}$, which holds for $\gamma\geq\frac{2N}{N+2}$. The strict inequality in (\ref{N-FMS}) results from the validity of the compact imbedding $\mathbf{V}_\gamma\hookrightarrow\hookrightarrow\mathbf{L}^2(\Omega)$, which is of fundamental use, not only in \cite{FMS-2003}, but also in all aforementioned works.

The analysis of the problem (\ref{geq1-inc})-(\ref{e-stress}) with $\mu_1=0$ and $\mu_2\not=0$ started, from the Mathematical Fluid Mechanics viewpoint, with the works on electrorheological fluids by R\.{u}\v{z}i\v{c}ka~\cite{R-1997} and on thermorheological fluids by Antontsev \emph{et al.}~\cite{ADO-PNDEA,AO-2011,AR-2006}. %
The main existence results for the problem (\ref{geq1-inc})-(\ref{e-stress}), with $\mu_1=0$ and $\mu_2\not=0$, are due to R\.{u}\v{z}i\v{c}ka~\cite{R-1997}, Huber~\cite{H-2011} and Diening \emph{et al.}~~\cite{DMS-2008}. %
The first existence result for the problem (\ref{geq1-inc})-(\ref{e-stress}), with $\mu_1=0$ and $\mu_2\not=0$, is an immediate consequence of the same result to the corresponding electrorheological problem established in~\cite{R-2000}. In fact, proceeding as in the proof of
\cite[Theorem 3.2.4]{R-2000}, one can easily proves the existence of weak solutions to our problem in the following class
\begin{equation}\label{Wq-nolog}
\mathbf{W}_{q(\cdot)}:=\mbox{closure of $\mathcal{V}$ in the $\|\mathbf{D(v)}\|_{\mathbf{L}^{q(\cdot)}(\Omega)}$--\ norm}\,,
\end{equation}
and satisfying to the energy relation (\ref{er-LL}) for all $\mathbf{\varphi}\in\mathbf{W}_{q(\cdot)}$. %
Moreover, the result holds for $q\in\mathcal{P}(\Omega)$ by assuming that
\begin{equation}\label{ap1}
\displaystyle
1<\alpha:=\mathrm{ess}\inf_{\hspace{-0.5cm}\mathbf{x}\in \Omega}q(\mathbf{x})\leq q(\mathbf{x})\leq \mathrm{ess}\sup_{\hspace{-0.5cm}\mathbf{x}\in \Omega}q(\mathbf{x}):=\beta<\infty,
\end{equation}
$\mathbf{f}\in\mathbf{V}_{\alpha}'$ and $\alpha$ satisfies to (\ref{N-LL}). %
The proof here follows the same approach of \cite{L-1967,Lions-1969} and uses the fact that $\mathbf{W}_{q(\cdot)}$ is imbedded into $\mathbf{V}_{\alpha}$. In~\cite{H-2011} the existence result~\cite{R-1997} was improved for the case of $\alpha$ satisfying to (\ref{N-FMS-R}). %
Here, under the assumptions that $\mathbf{f}\in(\mathbf{W}^{1,q(\cdot)}_0(\Omega))'$ and (\ref{glob-log-H}) holds, is proved the existence of weak solutions in the class
\begin{equation}\label{Vq}
\mathbf{V}_{q(\cdot)}:=\mbox{closure of $\mathcal{V}$ in $\mathbf{W}^{1,q(\cdot)}(\Omega)$}
\end{equation}
satisfying to the energy relation (\ref{er-LL}) for all $\mathbf{\varphi}\in\mathcal{V}$. %
The proof there combines a generalization of Bogowski\u{i}~\cite{B-1979} results on divergence problems in Orlicz-Sobolev spaces with the approach followed in \cite{FMS-1997,R-1997} for the case of a constant $q$. %
Finally, under the same assumptions of \cite{H-2011}, it is proved in~\cite{DMS-2008} a more complete existence result in the class $\mathbf{V}_{q(\cdot)}$. This result holds for $\alpha$ satisfying to (\ref{N-FMS})  and the solutions satisfies to the energy relation
\begin{equation}\label{er-DMS}
\int_{\Omega}\left(\mathbf{S}(\mathbf{D}(\mathbf{u}))-\mathbf{u}\otimes\mathbf{u}\right):\mathbf{D}(\mathbf{\varphi})\,d\mathbf{x}
=\int_{\Omega}p\,\mathrm{div}\varphi\,d\,\mathbf{x}+\int_{\Omega}\mathbf{f}\cdot\mathbf{\varphi}\,d\,\mathbf{x}
\end{equation}
for all $\mathbf{\varphi}\in\mathbf{W}_0^{1,\infty}(\Omega)$. %
The proof follows the same approach of the result for constant $q$ by the same authors and uses results on Lipschitz truncations of functions in Orlicz-Sobolev spaces performed by the authors still in~\cite{DMS-2008}.

A crucial assumption in the works~\cite{DMS-2008,H-2011} is that the variable exponent $q$ must be globally log-Hölder continuous. %
Recall that $q$ is globally log-Hölder continuous, if $q$ is locally log-Hölder continuous (\emph{f.} (\ref{glob-log-H})$_1$) and if there exist constants $C_2$ and $q_{\infty}$ such that (\ref{glob-log-H})$_2$ holds. %
This assumption gets into the proof in order to use the denseness of $\mathrm{C}_0^{\infty}(\Omega)$ in  $\mathrm{W}^{1,q(\cdot)}(\Omega)$ and, due to that, it is possible to look for the solutions in the class $\mathbf{V}_{q(\cdot)}$ instead of $\mathbf{W}_{q(\cdot)}$. %
However, though for many physical problems of trembling fluids, the variable power-law index $q$ satisfies the log-Hölder continuity property (\ref{glob-log-H}), there are some mathematical studies that suggest the non sufficiency of the condition (\ref{glob-log-H}). For a discussion on this important issue, see Diening \emph{et al.}~\cite[Section 5.1]{DHHR} and the references cited therein.
In the current work we shall establish an existence result for the complete problem (\ref{geq1-inc})-(\ref{e-stress}) without requiring the variable exponent $q$ is globally log-Hölder continuous.

\section{Weak formulation}\label{Sect-WF}

In order to introduce the notion of weak solutions we shall consider in this work, let us recall the well-known function spaces of Mathematical Fluid Mechanics defined at (\ref{V-reg})-(\ref{V-q=const}). Due to the presence of the variable exponent $q(\cdot)$ in the structure of the deviatoric tensor $\mathbf{S}$, we need to consider the weak solutions to the problem (\ref{geq1-inc})-(\ref{e-stress}) in some Orlicz-Sobolev space. %
As we already pointed out at Section~\ref{S-Prel}, an important limitation of Orlicz-Sobolev spaces is that, without any extra condition on the variable exponent $q$, but $q\in\mathcal{P}(\Omega)$ satisfying to (\ref{ap1}), the set $\mathbf{C}_0^{\infty}(\Omega)$ is not necessarily dense in $\mathbf{W}^{1,q(\cdot)}_0(\Omega)$. %
For this reason, we shall look for our solutions in the function space $\mathbf{W}_{q(\cdot)}$ defined in (\ref{Wq-nolog}) instead of the one considered in (\ref{Vq}). %
It is a easy task to verify the space $\mathbf{W}_{q(\cdot)}$ satisfies to the following imbeddings:
\begin{equation}\label{imb-a-b}
\mathbf{V}_{\beta}\hookrightarrow\mathbf{W}_{q(\cdot)}\hookrightarrow\mathbf{V}_{\alpha}\,.
\end{equation}
Moreover, $\mathbf{W}_{q(\cdot)}$ is a closed subspace of $\mathbf{V}_{\alpha}$ and therefore it is
a reflexive and separable Banach space for the norm
\begin{equation*}
\|\mathbf{v}\|_{\mathbf{W}_{q(\cdot)}}:=
\|\mathbf{D(v)}\|_{\mathbf{L}^{q(\cdot)}(\Omega)}.
\end{equation*}

\begin{definition}\label{weak-sol-vq}
Let $\Omega$ be a bounded domain of $\mathds{R}^N$, with $N \geq 2$. %
Assume that $\mathbf{f}\in\mathbf{L}^{1}(\Omega)$, $\gamma$ is a constant such that $1<\gamma<\infty$ and  $q\in\mathcal{P}(\Omega)$ is a variable exponent satisfying to (\ref{ap1}). %
A vector field $\mathbf{u}$ is a (very) weak solution to the problem (\ref{geq1-inc})-(\ref{e-stress}), if:
\begin{enumerate}
\item $\mathbf{\mathbf{u}}\in  \mathbf{W}_{q(\cdot)}\cap\mathbf{V}_{\gamma}$;
\item For every $\varphi\in\mathbf{W}_{q(\cdot)}\cap\mathbf{V}_{\gamma}$ (For every $\varphi\in\mathcal{V}$)
\begin{equation*}
\int_{\Omega}\left(\mu_1|\mathbf{D}(\mathbf{u})|^{\gamma-2}+\mu_2|\mathbf{D}(\mathbf{u})|^{q(\mathbf{x})-2}-\mathbf{u}\otimes\mathbf{u}\right):\mathbf{D}(\mathbf{\varphi})\,d\mathbf{x}
=\int_{\Omega}\mathbf{f}\cdot\varphi\,d\,\mathbf{x}.
\end{equation*}
\end{enumerate}
\end{definition}

\begin{remark}
Note that if $\gamma\geq\beta$, then $\mathbf{V}_{\gamma}\hookrightarrow\mathbf{W}_{q(\cdot)}$ and therefore it is enough to look for weak solutions in the class $\mathbf{W}_{q(\cdot)}$.
\end{remark}

As we had mentioned by the end of last section, in this work we are mainly focused on existence results for the problem (\ref{geq1-inc})-(\ref{e-stress}), or for some of its simplifications, without invoke the log-Hölder continuity property
(\ref{glob-log-H}) on the variable exponent $q$. We this in mind, an existence result for the problem that can be adapted from already known results is written in the following theorem.

\begin{theorem}\label{th0-exst-ws-pp}
Let $\Omega$ be a bounded domain in $\mathds{R}^{N}$, $N\geq 2$, with a Lipschitz-continuous boundary $\partial\Omega$.
Assume that $1<\gamma<\infty$, $q\in\mathcal{P}(\Omega)$ satisfies to (\ref{ap1}) and
$\mathbf{f}\in\left(\mathbf{V}_{\gamma}\cap\mathbf{W}_{q(\cdot)}\right)'$.
Then, if
\begin{equation*}
\min\left\{\gamma,\alpha\right\}\geq\frac{3N}{N+2}\,,
\end{equation*}
there exists a weak solution to the problem (\ref{geq1-inc})-(\ref{e-stress}) in the sense of Definition~\ref{weak-sol-vq}.
\end{theorem}
\begin{proof}
The proof combines the result addressed in \cite[Remark 2.5.5]{Lions-1969} for the constant power-law index $\gamma$ with the existence result of \cite[Theorem 3.2.4]{R-2000} for the variable power-law index $q$. We observe that, since (\ref{imb-a-b}) holds, \cite[Theorem 3.2.4]{R-2000} is a simple extension of the result established in \cite{Lions-1969}. Therefore, we can say that the proof Theorem~\ref{th0-exst-ws-pp} follows easily from the result invoked in \cite{Lions-1969}.
\end{proof}

\begin{remark}\label{R-Th0}
In the particular case of $\gamma=2$, Theorem~\ref{th0-exst-ws-pp} extends the existence result established by Ladyzhenskaya in~\cite{L-1967} to the case of a variable exponent $q$, without invoking the log-Hölder continuity property (\ref{glob-log-H}). %
Moreover, since $\mathbf{V}_{q(\cdot)}\subsetneqq\mathbf{W}_{q(\cdot)}$, this result is obtained in a larger class than the ones which follow from the works~\cite{DMS-2008,H-2011}.
\end{remark}

The main result of this work is written in what follows. We establish here the existence of weak solutions for the problem (\ref{geq1-inc})-(\ref{e-stress}) with a variable $q$ depending on $\mathbf{x}$ and without any other restriction on $q$ but
(\ref{ap1}). Our existence result will be valid for the constant exponent $\gamma$ and the variable one related via the relation (\ref{al-TE-1}) below.

\begin{theorem}\label{th-exst-ws-pp}
Let $\Omega$ be a bounded domain in $\mathds{R}^{N}$, $N\geq 2$. %
Assume that $q\in\mathcal{P}(\Omega)$ satisfies to (\ref{ap1}) and $\mathbf{f}\in\left(\mathbf{V}_{\gamma}\cap\mathbf{W}_{q(\cdot)}\right)'$.
Then, if for any $\delta>0$
\begin{equation}\label{al-TE-1}
\gamma\geq\max\left\{\frac{2N}{N+2}+\delta,\beta\right\}\,,
\end{equation}
there exists a very weak solution to the problem (\ref{geq1-inc})-(\ref{e-stress}) in the sense of Definition~\ref{weak-sol-vq}.
\end{theorem}

The proof of Theorem~\ref{th-exst-ws-pp} will be the aim of the next sections.  In order to make the proof as transparent as possible, we shall assume that
\begin{equation}\label{h-f}
\mathbf{f}=-\mathbf{div}\,\mathbf{F},\quad \mathbf{F}\in\mathds{M}^N_{\mathrm{sym}},
\quad \mathbf{F}\in\mathbf{L}^{q'(\cdot)}(\Omega),
\end{equation}
where $\mathds{M}^N_{\mathrm{sym}}$ is the vector space of all symmetric $N\times N$ matrices, which is equipped with the scalar product $\mathbf{A}:\mathbf{B}$ and
norm $|\mathbf{A}|=\sqrt{\mathbf{A}:\mathbf{A}}$.
Note that assumption (\ref{h-f}) does not affect at all the extent of Theorem~\ref{th-exst-ws-pp}, because $\mathbf{f}=-\mathbf{div}\,\mathbf{F}$ and $\mathbf{F}\in\mathbf{L}^{q'(\cdot)}(\Omega)$ implies that $\mathbf{f}\in\mathbf{W}_{q(\cdot)}'$, and $\gamma\geq\beta$ implies $\mathbf{W}_{q(\cdot)}'\hookrightarrow\left(\mathbf{V}_{\gamma}\cap\mathbf{W}_{q(\cdot)}\right)'$. %
The assumption $\mathbf{F}\in\mathds{M}^N_{\mathrm{sym}}$ is made in order to avoid unnecessary calculus. %
But, before we get into the details of the proof, let us make a few comments. %
Firstly, we note that, contrary to Theorem~\ref{th0-exst-ws-pp}, here we do not assume any regularity on the boundary $\partial\Omega$, besides the one following from the boundedness of $\Omega$.
Then, we observe that, as in Theorem~\ref{th0-exst-ws-pp}, the only assumption on $q$, besides (\ref{ap1}) and (\ref{al-TE-1}), is that $q$ is a measurable function from $\Omega$ into $[1,\infty]$.
Moreover, Remark~\ref{R-Th0} extends to this case of very weak solutions, now with $\gamma=2$ and $q$ related via (\ref{al-TE-1}).
This result is not mentioned in the literature, because, when $q$ is constant, Ladyzhenskaya's problem is, from the mathematical viewpoint, equivalent to the problem treated in Lions~\cite{Lions-1969}. This equivalence still remains in the case of a variable exponent $q$ satisfying to the log-Hölder continuity property (\ref{glob-log-H}) (\emph{f.} Diening \emph{et al.}~\cite{DMS-2008}). But if this hypothesis is not required, this result is, to the best of our knowledge, new. %
We observe that also the existence result established in Theorem~\ref{th-exst-ws-pp} is valid for a larger class of functions than the results provided by the works \cite{DMS-2008,H-2011}. %

\section{The regularized problem}\label{Sect-Exist-RP}

Let $\Phi\in\mathrm{C}^{\infty}([0,\infty))$ be a non-increasing function such that $0\leq\Phi\leq 1$ in $[0,\infty)$, $\Phi\equiv 1$ in $[0,1]$, $\Phi\equiv 0$ in $[2,\infty)$ and $0\leq -\Phi'\leq 2$. %
For $\epsilon>0$, we set
\begin{equation}\label{Phi-e}
\Phi_{\epsilon}(s):=\Phi(\epsilon s),\quad s\in[0,\infty).
\end{equation}
We consider the following regularized problem:
\begin{equation}\label{eq2-inc-e}
\mathrm{div}\,\mathbf{u}_{\epsilon}=0\quad\mbox{in}\quad \Omega,
\end{equation}
\begin{equation}\label{eq2-vel-e}
\begin{split}
\mathbf{div}(\mathbf{u}_{\epsilon}\otimes\mathbf{u}_{\epsilon}\Phi_{\epsilon}(|\mathbf{u}_{\epsilon}|))
=&\mathbf{f}-\mathbf{\nabla}p_{\epsilon}+\\
 &\mathbf{div}\left[\left(\mu_1|\mathbf{D}(\mathbf{u}_{\epsilon})|^{\gamma-2} +\mu_2|\mathbf{D}(\mathbf{u}_{\epsilon})|^{q(\mathbf{x})-2}\right)\mathbf{D}(\mathbf{u}_{\epsilon})\right]
\end{split}\quad\mbox{in}\quad \Omega,
\end{equation}
\begin{equation}\label{eq1-bc-u-e}
\mathbf{u}_{\epsilon}=\mathbf{0}\qquad\mbox{on}\quad \partial\Omega.
\end{equation}
Note that we have introduced a regularization that allows us to control the convective term. %
A vector function $\mathbf{u}_{\epsilon}\in\mathbf{V}_{\gamma}$ is a weak solution to the problem (\ref{eq2-inc-e})-(\ref{eq1-bc-u-e}), if
\begin{equation}\label{eq-ws-reg}
\begin{split}
&\int_{\Omega}\left[\left(\mu_1|\mathbf{D}(\mathbf{u}_{\epsilon})|^{\gamma-2}+
\mu_2|\mathbf{D}(\mathbf{u}_{\epsilon})|^{q(\mathbf{x})-2}\right)\mathbf{D}(\mathbf{u}_{\epsilon})-\mathbf{u}_{\epsilon}\otimes\mathbf{u}_{\epsilon}\Phi_{\epsilon}(|\mathbf{u}_{\epsilon}|)\right]:\mathbf{D}(\mathbf{\varphi})\,d\mathbf{x}\\
&=\int_{\Omega}\mathbf{F}:\mathbf{D}(\varphi)\,d\,\mathbf{x}
\end{split}
\end{equation}
for all $\varphi\in\mathcal{V}$. %
Note that (\ref{ap1}) implies $\mathbf{V}_{\gamma}\hookrightarrow\mathbf{W}_{q(\cdot)}$ and $\mathbf{f}$ is given in the form (\ref{h-f}).

\begin{proposition}\label{th-uxu-reg}
Let the assumptions of Theorem~\ref{th-exst-ws-pp} be fulfilled. %
Then, for each $\epsilon>0$, there exists a weak solution $\mathbf{u}_{\epsilon}\in\mathbf{V}_{\gamma}$ to the problem (\ref{eq2-inc-e})-(\ref{eq1-bc-u-e}). %
In addition, every weak solution satisfies to the following energy equality:
\begin{equation}\label{e-equality-qr}
\int_{\Omega}\left(\mu_1|\mathbf{D}(\mathbf{u}_{\epsilon})|^{\gamma}+
\mu_2|\mathbf{D}(\mathbf{u}_{\epsilon})|^{q(\mathbf{x})}\right)d\mathbf{x}
=\int_{\Omega}\mathbf{F}:\mathbf{D}(\mathbf{u}_{\epsilon})d\mathbf{x}.
\end{equation}
\end{proposition}

\begin{proof}
Let us set
$$M:=\{\mathbf{\varpi}\in\mathbf{H}:\|\mathbf{\varpi}\|_{\mathbf{L}^2(\Omega)}\leq 1\}$$
and lets us consider the following system:
\begin{equation}\label{eq2-inc-ve}
\mathrm{div}\,\mathbf{\upsilon}=0\quad\mbox{in}\quad \Omega,
\end{equation}
\begin{equation}\label{eq2-vel-ve}
\begin{split}
-\mathbf{div}\left[\left(\mu_1|\mathbf{D}(\mathbf{\upsilon})|^{\gamma-2}+\mu_2|\mathbf{D}(\mathbf{\upsilon})|^{q(\mathbf{x})-2}\right)\mathbf{D}(\mathbf{\upsilon})\right]
=&\mathbf{f}-\mathbf{\nabla}p\\
&-\mathbf{div}(\mathbf{\varpi}\otimes\mathbf{\varpi}\Phi_{\epsilon}(|\mathbf{\varpi}|))
\end{split}
\quad\mbox{in}\quad \Omega,
\end{equation}
\begin{equation}\label{eq1-bc-u-ve}
\mathbf{\upsilon}=\mathbf{0}\qquad\mbox{on}\quad \partial\Omega.
\end{equation}
Observing that, due to (\ref{al-TE-1}), the continuous imbedding $\mathbf{L}^{\gamma}(\Omega)\hookrightarrow\mathbf{L}^{q(\cdot)}(\Omega)$ holds, we can use the theory of monotone operators (see \emph{e.g.} Lions~\cite[Section 2.2]{Lions-1969}) to prove that for each $\mathbf{\varpi}\in M$, there exists a unique weak solution $\mathbf{\upsilon}\in\mathbf{V}_{\gamma}$ to the system (\ref{eq2-inc-ve})-(\ref{eq1-bc-u-ve}).

As a consequence of the previous step, we can define a mapping
\begin{equation}
\mathbf{K}:M\to\mathbf{H}
\end{equation}
such that to each $\mathbf{\varpi}\in M$ associates a unique $\mathbf{\varsigma}=\xi\,\mathbf{\upsilon}$, where $\mathbf{\upsilon}\in\mathbf{V}_{\gamma}$ is the unique weak solution to the system (\ref{eq2-inc-ve})-(\ref{eq1-bc-u-ve}) and $\xi$ is a positive constant which will be defined later on. %
Testing formally (\ref{eq2-vel-ve}) by the unique weak solution $\mathbf{\upsilon}$ such that $\xi\,\mathbf{\upsilon}:=\mathbf{K}(\mathbf{\varpi})$, with $\mathbf{\varpi}\in M$, integrating over $\Omega$, using the inequalities of Young and Korn, and at last the definition of $\Phi_{\epsilon}(|\mathbf{\varpi}|)$, we achieve to
\begin{equation}\label{est1-formaly}
\begin{split}
&C_1\int_{\Omega}\left(|\mathbf{D}(\mathbf{\upsilon})|^{\gamma}+
|\mathbf{D}(\mathbf{\upsilon})|^{q(\mathbf{x})}\right)d\mathbf{x}\leq\gamma_1+\gamma_2\|\mathbf{\varpi}\|_{\mathbf{L}^{2}(\Omega)}^2,
\\
&
\gamma_1:=C_2\int_{\Omega}|\mathbf{F}|^{\gamma'}d\mathbf{x},\quad
    \gamma_2:=C_3.
\end{split}
\end{equation}
Then, setting $\xi:=(C_4/C_1(\gamma_1+\gamma_2))^{-1}$, we can prove, using the inequalities of Sobolev and Korn, the generalized Hölder's inequality, and then (\ref{est1-formaly}), that
\begin{equation}\label{est-Kw-L2}
\|\mathbf{K}(\mathbf{\varpi})\|_{\mathbf{L}^2(\Omega)}\leq \xi\,C_4\left(\int_{\Omega}|\mathbf{D}(\mathbf{\upsilon})|^{q(\mathbf{x})}d\mathbf{x}+1\right)\leq \xi\,C_4/C_1(\gamma_1+\gamma_2)= 1
\end{equation}
for all $\mathbf{\varpi}\in M$. This proves that $\mathbf{K}$ maps $M$ into itself.

In order to prove the compactness of $\mathbf{K}$, we observe that from (\ref{est1-formaly}) it follows
\begin{equation}\label{est-Kw-Lg}
\|\mathbf{K}(\mathbf{\varpi})\|_{\mathbf{V}_{\gamma}}^{\gamma}\leq
C(\gamma_1+\gamma_2)
\end{equation}
for all $\mathbf{\varpi}\in M$. %
Owing to the embedding $\mathbf{L}^{q'(\cdot)}(\Omega)\hookrightarrow\mathbf{L}^{\gamma'}(\Omega)$ and to the assumption (\ref{h-f}), the right hand side of (\ref{est-Kw-Lg}) is finite.
Then, due to the compact imbedding $\mathbf{V}_{\gamma}\hookrightarrow\hookrightarrow\mathbf{L}^{2}(\Omega)$, valid for any $\gamma>\frac{2N}{N+2}$,  $\mathbf{K}(M)$ is relatively compact in $\mathbf{L}^{2}(\Omega)$.

To prove the continuity of $\mathbf{K}$, we consider a sequence $\mathbf{\varpi}_m$ in $M$ such that
$$\mathbf{\varpi}_m\to \mathbf{\varpi}\quad\mbox{in}\quad \mathbf{L}^{2}(\Omega)\quad\mbox{as}\quad m\to\infty.$$
By the relative compactness of $\mathbf{K}(M)$ in $\mathbf{L}^{2}(\Omega)$, there exists a subsequence $\mathbf{\varpi}_{m_k}$ such that
\begin{equation}\label{cont-subs-reg}
\mathbf{K}(\mathbf{\varpi}_{m_k})\to \mathbf{\varsigma}\quad\mbox{in}\quad \mathbf{L}^{2}(\Omega),\quad\mbox{as}\quad k\to\infty.
\end{equation}
From the definition of $\mathbf{K}$, the functions $\mathbf{\upsilon}_{m_k}$ defined by
$\xi\,\mathbf{\upsilon}_{m_k}\equiv\varsigma_{m_k}:=\mathbf{K}(\mathbf{\varpi}_{m_k})$, satisfy to
\begin{equation}\label{eq-ws-seq-reg}
\begin{split}
  & \int_{\Omega}\left(\mu_1|\mathbf{D}(\mathbf{\upsilon}_{m_k})|^{\gamma-2}+
   \mu_2|\mathbf{D}(\mathbf{\upsilon}_{m_k})|^{q(\mathbf{x})-2}\right)\mathbf{D}(\mathbf{\upsilon}_{m_k}):\mathbf{D}(\mathbf{\varphi})\,d\mathbf{x}
= \\
  & \int_{\Omega}\left(\mathbf{F}+\mathbf{\varpi}_{m_k}\otimes\mathbf{\varpi}_{m_k}\Phi_{\epsilon}(|\mathbf{\varpi}_{m_k}|)\right):\mathbf{D}(\varphi)\,d\,\mathbf{x}
\end{split}
\end{equation}
for all $\varphi\in\mathcal{V}$. %
Passing to the limit in (\ref{eq-ws-seq-reg}) and using the Minty trick, we can prove that
$\varsigma=\xi\,\mathbf{\upsilon}$ and therefore $\varsigma=\mathbf{K}(\mathbf{\varpi})$. %
Then, from (\ref{cont-subs-reg}), we conclude that
$\mathbf{\mathbf{K}}(\mathbf{\varpi}_m)\to \mathbf{K}(\mathbf{\varpi})$ in $\mathbf{L}^{2}(\Omega)$ as $m\to\infty$, which proves the continuity of $\mathbf{K}$.

Now, applying Schauder's fixed point theorem, there exists a function $\mathbf{\upsilon}_{\xi}\in M$  such that $\mathbf{K}(\mathbf{\upsilon}_{\xi})=\mathbf{\upsilon}_{\xi}$ and which is a weak solution to the problem (\ref{eq2-inc-e})-(\ref{eq1-bc-u-e}) in $\Omega$.

Finally, the energy relation (\ref{e-equality-qr}) follows by testing (\ref{eq2-vel-e}) by a weak solution and integrating over $\Omega$, and observing that now the convective term is zero.
\end{proof}

\section{Existence of approximative solutions}\label{Sect-Exist-AS}

Let $\mathbf{u}_{\epsilon}\in\mathbf{V}_{\gamma}$ be a weak solution to the problem (\ref{eq2-inc-e})-(\ref{eq1-bc-u-e}). %
From Proposition~\ref{th-uxu-reg} (\emph{f.} (\ref{e-equality-qr})), using the generalized Hölder's inequality and the relation between the norm $\|\mathbf{D}(\mathbf{u}_{\epsilon})\|_{\mathbf{L}^{q(\cdot)}(\Omega)}$ and the semimodular $A_{q(\cdot)}(\mathbf{D}(\mathbf{u}_{\epsilon}))$ (\emph{f.} (\ref{ap4-56})) together with Young's inequality, we can prove that
\begin{equation}\label{e-inequality-qr}
\int_{\Omega}\left(|\mathbf{D}(\mathbf{u}_{\epsilon})|^{\gamma}+
|\mathbf{D}(\mathbf{u}_{\epsilon})|^{q(\mathbf{x})}\right)d\mathbf{x}\leq C,
\end{equation}
where, by the assumption (\ref{h-f}), $C$ is a positive constant and, very important, does not depend on $\epsilon$. %
Appealing again to the relation between the norm $\|\mathbf{D}(\mathbf{u}_{\epsilon})\|_{\mathbf{L}^{q(\cdot)}(\Omega)}$ and the semimodular $A_{q(\cdot)}(\mathbf{D}(\mathbf{u}_{\epsilon}))$ and also to generalized Hölder's inequality, we can prove from (\ref{e-inequality-qr}) that
\begin{equation}\label{est-inf-gam}
\|\mathbf{u}_{\epsilon}\|_{\mathbf{V}_{\gamma}}\leq C,
\end{equation}
\begin{equation}\label{est-q}
\|\mathbf{D}(\mathbf{u}_{\epsilon})\|_{\mathbf{L}^{q(\cdot)}(\Omega)}\leq C.
\end{equation}
Proceeding as for (\ref{est-inf-gam})-(\ref{est-q}), we can also prove that
\begin{equation}\label{est-gam'}
\||\mathbf{D}(\mathbf{u}_{\epsilon})|^{\gamma-2}\mathbf{D}(\mathbf{u}_{\epsilon})\|_{\mathbf{L}^{\gamma'}(\Omega)}\leq C,
\end{equation}
\begin{equation}\label{est-q'}
\||\mathbf{D}(\mathbf{u}_{\epsilon})|^{q(\mathbf{x})-2}\mathbf{D}(\mathbf{u}_{\epsilon})\|_{\mathbf{L}^{q'(\cdot)}(\Omega)}\leq C.
\end{equation}
Moreover, using (\ref{al-TE-1}), it also follows that
\begin{equation}\label{est-q'g'}
\||\mathbf{D}(\mathbf{u}_{\epsilon})|^{q(\mathbf{x})-2}\mathbf{D}(\mathbf{u}_{\epsilon})\|_{\mathbf{L}^{\gamma'}(\Omega)}\leq C.
\end{equation}
On the other hand, by using (\ref{est-inf-gam}) and Sobolev's inequality, we have
\begin{equation}\label{est-gam(N+2)/N}
\|\mathbf{u}_{\epsilon}\|_{\mathbf{L}^{\gamma^{\ast}}(\Omega)}\leq C,
\end{equation}
where $\gamma^{\ast}$ denotes the Sobolev conjugate of $\gamma$. %
As a consequence of (\ref{est-gam(N+2)/N}) and due to the definition of $\Phi_{\epsilon}$ (\emph{f.} (\ref{Phi-e})),
\begin{equation}\label{est-Phiuxu}
\|\mathbf{u}_{\epsilon}\otimes\mathbf{u}_{\epsilon}\Phi_{\epsilon}(|\mathbf{u}_{\epsilon}|)\|_{\mathbf{L}^{\frac{\gamma^{\ast}}{2}}(\Omega)}\leq C.
\end{equation}
Note that the constants in (\ref{est-inf-gam})-(\ref{est-Phiuxu}) are distinct and do not depend on $\epsilon$. %
From (\ref{est-inf-gam}), (\ref{est-gam'}), (\ref{est-q'g'}) and (\ref{est-Phiuxu}), there exists a sequence of positive numbers $\epsilon_m$ such that
$\epsilon_m\to 0$, as $m\to \infty$, and
%\begin{equation}\label{conv-Lq}
%\mathbf{D}(\mathbf{u}_{\epsilon_m})\to\mathbf{D}(\mathbf{u})\quad\mbox{weakly in}\quad\mathbf{L}^{q(\cdot)}(\Omega),\quad\mbox{as}\ m\to\infty,
%\end{equation}
%\begin{equation}\label{convg-S-q'}
%\mbox{$\mathbf{S}(\mathbf{D}(\mathbf{u}_{\epsilon_m}))\to \mathbf{S}$\quad weakly in $\mathbf{L}^{q'(\cdot)}(\Omega)$,\quad as $m\to\infty$,}
%\end{equation}
\begin{equation}\label{convg-La}
\mbox{$\mathbf{u}_{\epsilon_m}\to \mathbf{u}$\quad weakly in $\mathbf{V}_{\gamma}$,\quad as $m\to\infty$,}
\end{equation}
\begin{equation}\label{conv-Re-0}
|\mathbf{D}(\mathbf{u}_{\epsilon_m})|^{\gamma-2}\mathbf{D}(\mathbf{u}_{\epsilon_m})\to \mathbf{S}_{1}\quad\mbox{weakly in}\quad\mathbf{L}^{\gamma'}(\Omega),\quad\mbox{as}\quad m\to\infty.
\end{equation}
\begin{equation}\label{convg-S-b'}
|\mathbf{D}(\mathbf{u}_{\epsilon_m})|^{q(\mathbf{x})-2}\mathbf{D}(\mathbf{u}_{\epsilon_m})\to \mathbf{S}_{2}\quad\mbox{weakly in}\quad\mathbf{L}^{\gamma'}(\Omega),\quad\mbox{as}\quad m\to\infty.
\end{equation}
%\begin{equation}\label{convg-La(N+2)/N}
%\mbox{$\mathbf{u}_{\epsilon_m}\to \mathbf{u}$\quad weakly in $\mathbf{L}^{\alpha^{\ast}}(\Omega)$,\quad as $m\to\infty$,}
%\end{equation}
\begin{equation}\label{convg-Phiuxu}
\mbox{$\mathbf{u}_{\epsilon_m}\otimes\mathbf{u}_{\epsilon_m}\Phi_{\epsilon_m}(|\mathbf{u}_{\epsilon_m}|)\to \mathbf{G}$\quad weakly in $\mathbf{L}^{\frac{\gamma^{\ast}}{2}}(\Omega)$,\quad as $m\to\infty$}.
\end{equation}

Now we observe that, due to (\ref{convg-La}), the application of Sobolev's compact imbedding theorem implies
\begin{equation}\label{convg-s-g}
\mbox{$\mathbf{u}_{\epsilon_m}\to \mathbf{u}$\quad strongly in $\mathbf{L}^{\kappa}(\Omega)$,\quad as $m\to\infty$,
\quad
for any $\kappa:1\leq\kappa<\gamma^{\ast}$.
}
\end{equation}
Since (\ref{al-TE-1}) implies $2<\gamma^{\ast}$, it follows from (\ref{convg-s-g}) that
\begin{equation}\label{convg-s-2}
\mbox{$\mathbf{u}_{\epsilon_m}\to \mathbf{u}$\quad strongly in $\mathbf{L}^{2}(\Omega)$,\quad as $m\to\infty$.
}
\end{equation}
Using the definition of $\Phi_{\epsilon_m}$ (\emph{f.} (\ref{Phi-e})) and the result (\ref{convg-s-2}), we can prove that
\begin{equation}\label{limit-eq-G}
\mbox{$\mathbf{u}_{\epsilon_m}\otimes\mathbf{u}_{\epsilon_m}\Phi_{\epsilon_m}(|\mathbf{u}_{\epsilon_m}|)\to \mathbf{u}\otimes\mathbf{u}$\quad strongly in $\mathbf{L}^{1}(\Omega)$,\quad as $m\to\infty$.
}
\end{equation}
Then gathering the information of (\ref{convg-Phiuxu}) and (\ref{limit-eq-G}), we see that $\mathbf{G}=\mathbf{u}\otimes\mathbf{u}$.

Finally, using the convergence results (\ref{convg-La})-(\ref{convg-Phiuxu}) and observing (\ref{limit-eq-G}), we can pass to the limit $m\to\infty$ in the following integral identity, which results from (\ref{eq-ws-reg}),
\begin{equation}\label{eq-wsm-reg}
\begin{split}
&\int_{\Omega}\left(\mu_1|\mathbf{D}(\mathbf{u}_{\epsilon_m})|^{\gamma-2}+\mu_2|\mathbf{D}(\mathbf{u}_{\epsilon_m})|^{q(\mathbf{x})-2}\right)\mathbf{D}(\mathbf{u}_{\epsilon_m}):\mathbf{D}(\mathbf{\varphi})\,d\mathbf{x}\\
&-\int_{\Omega}\left[\mathbf{u}_{\epsilon_m}\otimes\mathbf{u}_{\epsilon_m}\Phi_{\epsilon_m}(|\mathbf{u}_{\epsilon_m}|)+\mathbf{F}\right]:\mathbf{D}(\mathbf{\varphi})\,d\mathbf{x}
=0,
\end{split}
\end{equation}
valid for all $\mathbf{\varphi}\in\mathcal{V}$, to obtain
\begin{equation}\label{limit-eq-SH}
\int_{\Omega}(\mu_1\mathbf{S}_{1}+\mu_2\mathbf{S}_{2}-\mathbf{u}\otimes\mathbf{u}-\mathbf{F}):\mathbf{D}(\mathbf{\varphi})\,d\mathbf{x}=0\quad
\forall\ \varphi\in\mathcal{V}. %
\end{equation}

\section{Determination of the pressure}\label{Sect-DP}

Since we shall use test functions which are not divergence free, we first have to determine the approximative pressure from the weak formulation
(\ref{eq-wsm-reg}). %
First, let $\omega'$ be a fixed but arbitrary open bounded subset of $\Omega$ such that
\begin{equation}\label{om-l}
\omega'\subset\subset\Omega\quad\mbox{and}\quad\partial\omega'\ \mbox{is Lipschitz} %
\end{equation}
and let us set
\begin{equation}\label{Q-em}
\mathbf{Q}_{\epsilon_m}:=\left(\mu_1|\mathbf{D}(\mathbf{u}_{\epsilon_m})|^{\gamma-2}+\mu_2|\mathbf{D}(\mathbf{u}_{\epsilon_m})|^{q(\mathbf{x})-2}\right)\mathbf{D}(\mathbf{u}_{\epsilon_m})-\mathbf{u}_{\epsilon_m}\otimes\mathbf{u}_{\epsilon_m}\Phi_{\epsilon_m}(|\mathbf{u}_{\epsilon_m}|)-\mathbf{F}.
\end{equation}
Using assumption (\ref{h-f}) and the results (\ref{est-gam'}), (\ref{est-q'g'}), (\ref{est-Phiuxu}) and (\ref{conv-Re-0}), we can prove that
\begin{equation}\label{Q-in-Lr}
\mathbf{Q}_{\epsilon_m}\in\mathbf{L}^r(\Omega),
\end{equation}
where $r$ can be taken in such a way that
\begin{equation}\label{r}
 1<r\leq r_0:=\min\left\{\gamma',\frac{\gamma^{\ast}}{2}\right\}.
\end{equation}
Note that $r_0=\min\{\gamma',\frac{N\gamma}{2(N-\gamma)}\}$ if $\gamma<N$ and we can take $r_0=\gamma'$ if $N\geq\gamma$. %
Next, we define a linear functional
\begin{equation}\label{b1-lin-f-r}
\Pi_{\epsilon_m}:\mathbf{W}^{1,r'}_0(\omega')\rightarrow\mathbf{W}^{-1,r}(\omega')
\end{equation}
by
\begin{equation}\label{b2-lin-f-r}
\langle\Pi_{\epsilon_m},\mathbf{\varphi}\rangle_{\mathbf{W}^{-1,r}(\omega')\times\mathbf{W}^{1,r'}_0(\omega')}:=
\int_{\omega'}\mathbf{Q}_{\epsilon_m}:\mathbf{D}(\mathbf{\varphi})\,d\mathbf{x}.
\end{equation}
Using (\ref{b1-lin-f-r})-(\ref{b2-lin-f-r}), we can prove, owing to (\ref{Q-in-Lr}), that exists a positive constant $C$ independent of $m$ such that
\begin{equation}\label{b-lin-op-p}
\|\Pi_{\epsilon_m}\|_{(\mathbf{V}_{r'})'}\leq C.
\end{equation}
Note that here $\mathbf{V}_{r'}$ is taken over $\omega'$. %
Moreover, since $\mathcal{V}$ is dense in $\mathbf{V}_{r'}$, we can see,  due to (\ref{eq-wsm-reg}), (\ref{Q-em}) and (\ref{b2-lin-f-r}), that
\begin{equation}\label{dual=0}
\langle\Pi_{\epsilon_m},\mathbf{\varphi}\rangle_{(\mathbf{V}_{r'})'\times\mathbf{V}_{r'}}=0\quad\forall\
\mathbf{\varphi}\in\mathbf{V}_{r'}.
\end{equation}
By virtue of (\ref{b1-lin-f-r})-(\ref{dual=0}) and due to assumption (\ref{om-l}), we can apply a version of de Rham's Theorem (\emph{f}. Bogovski\u{i}~\cite[Theorems 1-4]{B-1980} and
Pileckas~\cite[Section 1]{P-1983}) to prove the existence of a unique function
\begin{equation}\label{app-pr}
p_{\epsilon_m}\in\mathbf{L}^{r'}(\omega'),\qquad\mbox{with}\quad\int_{\omega'}p_{\epsilon_m}d\mathbf{x}=0,
\end{equation}
such that
\begin{equation}\label{Rham-thm}
\langle\Pi_{\epsilon_m},\mathbf{\varphi}\rangle_{\mathbf{W}^{-1,r}(\omega')\times\mathbf{W}^{1,r'}_0(\omega')}=
\int_{\omega'}p_{\epsilon_m}\,\mathrm{div}\mathbf{\varphi}\,d\mathbf{x}\quad\forall\ \mathbf{\varphi}\in\mathbf{W}^{1,r'}_0(\omega')
\end{equation}
and
\begin{equation}\label{b-pem}
\|p_{\epsilon_m}\|_{\mathbf{L}^{r'}(\omega')}\leq\|\Pi_{\epsilon_m}\|_{(\mathbf{V}_{r'})'}.
\end{equation}
Then, gathering the information of (\ref{eq-wsm-reg}), (\ref{Q-em}), (\ref{b2-lin-f-r}) and (\ref{Rham-thm}), we obtain
\begin{equation}\label{eq-wsm-reg-pem}
\begin{split}
 & \int_{\omega'}\left(\mu_1|\mathbf{D}(\mathbf{u}_{\epsilon_m})|^{\gamma-2}+
\mu_2|\mathbf{D}(\mathbf{u}_{\epsilon_m})|^{q(\mathbf{x})-2}\right)\mathbf{D}(\mathbf{u}_{\epsilon_m}):\mathbf{D}(\mathbf{\varphi})\,d\mathbf{x}= \\
 & \int_{\omega'}\mathbf{F}:\mathbf{D}(\mathbf{\varphi})\,d\mathbf{x}+
\int_{\omega'}\mathbf{u}_{\epsilon_m}\otimes\mathbf{u}_{\epsilon_m}\Phi_{\epsilon_m}(|\mathbf{u}_{\epsilon_m}|):\mathbf{D}(\mathbf{\varphi})\,d\mathbf{x}+
\int_{\omega'}p_{\epsilon_m}\,\mathrm{div}\mathbf{\varphi}\,d\mathbf{x}
\end{split}
\end{equation}
for all $\mathbf{\varphi}\in\mathbf{W}^{1,r'}_0(\omega')$. %
On the other hand, due to (\ref{b-lin-op-p}) and (\ref{b-pem}) and by means of reflexivity, we get, passing to a subsequence, that
\begin{equation}\label{conv-p0-em}
p_{\epsilon_m}\to p_0\quad\mbox{weakly in}\quad\mathrm{L}^{r'}(\omega'),\quad\mbox{as}\quad m\to\infty.
\end{equation}
Next, passing to the limit $m\to\infty$ in the integral identity (\ref{eq-wsm-reg-pem}) by using the convergence results (\ref{conv-Re-0})-(\ref{convg-Phiuxu}), observing that in the last case, by virtue of (\ref{limit-eq-G}), $\mathbf{G}=\mathbf{u}\otimes\mathbf{u}$, using also (\ref{conv-p0-em}), we obtain
\begin{equation}\label{eq-limit-pem}
\begin{split}
 & \int_{\omega'}\left(\mu_1\mathbf{S}_1+\mu_2\mathbf{S}_2-\mathbf{u}\otimes\mathbf{u}-\mathbf{F}\right):\mathbf{D}(\mathbf{\varphi})\,d\mathbf{x}=
\int_{\omega'}p_0\,\mathrm{div}\mathbf{\varphi}\,d\mathbf{x}
\end{split}
\end{equation}
for all $\mathbf{\varphi}\in\mathbf{W}^{1,r'}_0(\omega')$. %
Then, proceeding analogously as we did for (\ref{b1-lin-f-r})-(\ref{b-lin-op-p}), we can define a linear functional
\begin{equation}\label{b1-lin-0}
\Pi_0:\mathbf{W}^{1,r'}_0(\omega')\rightarrow\mathbf{W}^{-1,r}(\omega')
\end{equation}
by
\begin{equation}\label{b2-lin-0}
\langle\Pi_0,\mathbf{\varphi}\rangle_{\mathbf{W}^{-1,r}(\omega')\times\mathbf{W}^{1,r'}_0(\omega')}:=
\int_{\omega'}\mathbf{Q}_0:\mathbf{D}(\mathbf{\varphi})\,d\mathbf{x},
\end{equation}
where $\mathbf{Q}_0:=\mu_1\mathbf{S}_1+\mu_2\mathbf{S}_2-\mathbf{u}\otimes\mathbf{u}-\mathbf{F}$, such that
(\ref{b-lin-op-p})-(\ref{dual=0}) are verified with $\Pi_0$ and $\mathbf{Q}_0$ in the places of $\Pi_{\epsilon_m}$ and $\mathbf{Q}_{\epsilon_m}$.
In consequence, by the same version of de Rham's Theorem aforementioned, there exists a unique function
\begin{equation}\label{app-pr-0}
\overline{p}_0\in\mathbf{L}^{r'}(\omega'),\qquad\mbox{with}\quad\int_{\omega'}\overline{p}_0\,d\mathbf{x}=0,
\end{equation}
such that
\begin{equation}\label{Rham-thm-0}
\langle\Pi_0,\mathbf{\varphi}\rangle_{\mathbf{W}^{-1,r}(\omega')\times\mathbf{W}^{1,r'}_0(\omega')}=
\int_{\omega'}\overline{p}_0\,\mathrm{div}\mathbf{\varphi}\,d\mathbf{x},\quad\forall\ \mathbf{\varphi}\in\mathbf{W}^{1,r'}_0(\omega')
\end{equation}
and (\ref{b-pem}) is verified with $\overline{p}_0$ and $\Pi_0$ in the places of $p_{\epsilon_m}$ and $\Pi_{\epsilon_m}$. %
Then gathering (\ref{b2-lin-0}) and (\ref{Rham-thm-0}), we achieve to
\begin{equation}\label{eq-wsm-reg-p0}
\int_{\omega'}\left(\mu_1\mathbf{S}_1+\mu_2\mathbf{S}_2-\mathbf{u}\otimes\mathbf{u}-\mathbf{F}\right):\mathbf{D}(\mathbf{\varphi})\,d\mathbf{x}=
\int_{\omega'}
\overline{p}_0\,\mathrm{div}\mathbf{\varphi}\,d\mathbf{x}
\end{equation}
for all $\mathbf{\varphi}\in\mathbf{W}^{1,r'}_0(\omega')$. %
Finally, combining (\ref{eq-limit-pem}) and (\ref{eq-wsm-reg-p0}), and by means of uniqueness, we conclude that
\begin{equation*}
\overline{p}_0=p_0.
\end{equation*}

\section{Decomposition of the pressure.}\label{Sect-Dec-p}

The main idea in this section is the application of a method to locally decompose the pressure found in the previous section. %
For that, we shall use a lemma which is proved by using a direct decomposition of $L^s$, which in turn is equivalent to the weak $L^s$-solvability
of the Dirichlet problem for the Bilaplacian in bounded domains with $C^2$ boundaries (\emph{f.}
Simader and Sohr~\cite{SS-1996}).
%, and by using a result on a variational estimate (\emph{f.} ~\cite{SS-1996}). COLOCAR AQUI AS REFERÊNCIAS EXACTAS DOS TEOREMAS DESTE LIVRO
With this in mind, let  $\omega$ be a fixed but arbitrary domain such that
\begin{equation}\label{om-C2}
\omega\subset\subset\omega'\subset\subset\Omega\quad\mbox{and}\quad\partial\omega\ \mbox{is $C^2$}. %
\end{equation}

\begin{lemma}\label{L-2.4-Wolf}
Let $1<s<\infty$. Then for every $v^{\ast}\in\mathrm{W}^{-2,s}(\omega)$ there exists a unique $u\in\mathrm{W}^{2,s}_0(\omega)$ such that
\begin{equation*}\label{dec-Lem}
\int_{\omega}\triangle u\,\triangle\phi\,d\mathbf{x}=
\langle v^{\ast},\phi\rangle_{\mathbf{W}^{-2,s}(\omega)\times\mathbf{W}^{2,s'}_0(\omega)}\qquad\forall\ \phi\in\mathrm{C}_0^{\infty}(\omega).
\end{equation*}
\end{lemma}
\begin{proof}
See \emph{e.g.} Wolf~\cite[Lemma 2.4]{Wolf-2007}.
\end{proof}
To simplify the notation in the sequel, let us set
\begin{equation*}\label{set-A}
\mathrm{A}^s(\omega):=\{a\in\mathrm{L}^{s}(\omega):a=\triangle u,\quad u\in\mathrm{W}^{2,s}_0(\omega)\},\qquad  1<s<\infty. \\
\end{equation*}
Applying Lemma~\ref{L-2.4-Wolf}, with $s=\gamma'$ first and then with $s=\frac{\gamma^{\ast}}{2}$, attending to the definitions of $\mathrm{A}^{\gamma'}(\omega)$ and $\mathrm{A}^{\frac{\gamma^{\ast}}{2}}(\omega)$, and using (\ref{conv-Re-0}) and (\ref{convg-S-b'}) by one hand and (\ref{convg-Phiuxu}) and (\ref{limit-eq-G}) on the other, we can infer that exist unique functions
\begin{equation}\label{E1-p1}
p^1_{\epsilon_m}\in\mathrm{A}^{\gamma'}(\omega),
\end{equation}
\begin{equation}\label{E1-p2}
p^2_{\epsilon_m}\in\mathrm{A}^{\frac{\gamma^{\ast}}{2}}(\omega)
\end{equation}
such that
\begin{equation}\label{dec-p1}
\begin{split}
&\int_{\omega}p^1_{\epsilon_m}\triangle\phi\,d\mathbf{x}=\\
&\int_{\omega}\left[\mu_1\left(|\mathbf{D}(\mathbf{u}_{\epsilon_m})|^{\gamma-2}\mathbf{D}(\mathbf{u}_{\epsilon_m})-\mathbf{S}_1\right)+\mu_2\left(|\mathbf{D}(\mathbf{u}_{\epsilon_m})|^{q(\mathbf{x})-2}\mathbf{D}(\mathbf{u}_{\epsilon_m})- \mathbf{S}_2\right)\right]:\nabla^2\phi\,d\mathbf{x},
\end{split}
\end{equation}
\begin{equation}\label{dec-p2}
\int_{\omega}p^2_{\epsilon_m}\triangle\phi\,d\mathbf{x}=-
\int_{\omega}(\mathbf{u}_{\epsilon_m}\otimes\mathbf{u}_{\epsilon_m}\Phi_{\epsilon_m}(|\mathbf{u}_{\epsilon_m}|)-
\mathbf{u}\otimes\mathbf{u}):\nabla^2\phi\,d\mathbf{x}
\end{equation}
for all $\phi\in\mathrm{C}_0^{\infty}(\omega)$. %
The next result helps us with the estimates of the local pressures found in (\ref{E1-p1})-(\ref{E1-p2}).
\begin{lemma}\label{L-2.3-Wolf}
Let $1<s<\infty$ and assume that $p\in\mathrm{L}^s(\omega)$ and $\mathbf{Q}\in\mathbf{L}^{s}(\omega)$ are such that
\begin{equation*}\label{dec-p}
\int_{\omega}p\triangle\phi\,d\mathbf{x}=
\int_{\omega}\mathbf{Q}:\nabla^2\phi\,d\mathbf{x}\quad\forall\ \phi\in\mathrm{C}^{\infty}_0(\omega).
\end{equation*}
Then
\begin{equation*}\label{b-p}
\|p\|_{\mathrm{L}^s(\omega)}\leq
C\|\mathbf{Q}\|_{\mathbf{L}^s(\omega)},
\end{equation*}
where $C$ is a positive constant depending on $s$, $N$, $\omega$ and on the Calderón-Zigmund inequality's constant.
\end{lemma}
\begin{proof}
See Wolf~\cite[Lemma 2.3]{Wolf-2007}.
\end{proof}

\noindent Attending to (\ref{convg-S-b'}), (\ref{conv-Re-0}) and (\ref{dec-p1}) by one hand, and (\ref{convg-Phiuxu}), (\ref{limit-eq-G}) and (\ref{dec-p2}) on the other, a direct application of Lemma~\ref{L-2.3-Wolf}, with $s=\gamma'$ and then with $s=\frac{\gamma^{\ast}}{2}$, yields
\begin{equation}\label{bd2-p1-em}
\begin{split}
\|p^1_{\epsilon_m}\|_{\mathrm{L}^{\gamma'}(\omega)}\leq & C_1\||\mathbf{D}(\mathbf{u}_{\epsilon_m})|^{\gamma-2}\mathbf{D}(\mathbf{u}_{\epsilon_m})-\mathbf{S}_1\|_{\mathbf{L}^{\gamma'}(\omega)}+\\
& C_2\||\mathbf{D}(\mathbf{u}_{\epsilon_m})|^{q(\mathbf{x})-2}\mathbf{D}(\mathbf{u}_{\epsilon_m})-\mathbf{S}_2\|_{\mathbf{L}^{\gamma'}(\omega)},
\end{split}
\end{equation}
\begin{equation}\label{bd2-p2-em}
\|p^2_{\epsilon_m}\|_{\mathrm{L}^{\frac{\gamma^{\ast}}{2}}(\omega)}\leq C_3\|\mathbf{u}_{\epsilon_m}\otimes\mathbf{u}_{\epsilon_m}\Phi_{\epsilon_m}(|\mathbf{u}_{\epsilon_m}|)-
\mathbf{u}\otimes\mathbf{u}\|_{\mathbf{L}^{\frac{\gamma^{\ast}}{2}}(\omega)}.
\end{equation}
where $C_1$, $C_2$ and $C_3$ are positive constants depending on $\gamma$, $N$, $\omega$ and on the Calderón-Zigmund inequality's constant. %

On the other hand, combining (\ref{eq-wsm-reg-pem}) and (\ref{eq-limit-pem}), and using the definition of the distributive derivative, we obtain
\begin{equation}\label{dist-der}
\begin{split}
   & \mathbf{div}\left(|\mathbf{D}(\mathbf{u}_{\epsilon_m})|^{\gamma-2}\mathbf{D}(\mathbf{u}_{\epsilon_m})-\mathbf{S}_1+|\mathbf{D}(\mathbf{u}_{\epsilon_m})|^{q(\mathbf{x})-2}\mathbf{D}(\mathbf{u}_{\epsilon_m})-\mathbf{S}_2\right)-\\
   &\mathbf{div}\left(\mathbf{u}_{\epsilon_m}\otimes\mathbf{u}_{\epsilon_m}\Phi_{\epsilon_m}(|\mathbf{u}_{\epsilon_m}|)-\mathbf{u}\otimes\mathbf{u}\right)
   =\mathbf{\nabla}(p_{\epsilon_m}-p_0)
\end{split}
\qquad\mbox{in $\mathcal{D}'(\omega)$}.
\end{equation}
Then, testing (\ref{dist-der}) by $\nabla\phi$, with $\phi\in\mathrm{C}_0^{\infty}(\omega)$, integrating over $\omega$ and comparing the resulting equation with the one resulting from adding (\ref{dec-p1}) and (\ref{dec-p2}), we obtain
$$p_{\epsilon_m}-p_0=p^1_{\epsilon_m}+p^2_{\epsilon_m}\,.$$
Inserting this into (\ref{dist-der}), it follows that
\begin{equation}\label{dist2-der}
\begin{split}
   & \mathbf{div}\left(|\mathbf{D}(\mathbf{u}_{\epsilon_m})|^{\gamma-2}\mathbf{D}(\mathbf{u}_{\epsilon_m})-\mathbf{S}_1+|\mathbf{D}(\mathbf{u}_{\epsilon_m})|^{q(\mathbf{x})-2}\mathbf{D}(\mathbf{u}_{\epsilon_m})-\mathbf{S}_2\right)-\\
   &\mathbf{div}\left(\mathbf{u}_{\epsilon_m}\otimes\mathbf{u}_{\epsilon_m}\Phi_{\epsilon_m}(|\mathbf{u}_{\epsilon_m}|)-\mathbf{u}\otimes\mathbf{u}\right)
   =\mathbf{\nabla}\left(p^1_{\epsilon_m}+p^2_{\epsilon_m}\right)
\end{split}
\qquad\mbox{in $\mathcal{D}'(\omega)$}.
\end{equation}

\section{The Lipschitz truncation}\label{Sect-Lip-T}

To start this section, let us set
\begin{equation}\label{v-em-chi}
\mathbf{w}_{\epsilon_m}:=(\mathbf{u}_{\epsilon_m}-\mathbf{u})\chi_{\omega},
\end{equation}
where $\chi_{\omega}$ denotes the characteristic function of the set $\omega$ introduced in (\ref{om-C2}). %
Having in mind the extension of (\ref{dist2-der}) to $\mathds{R}^{N}$, here we shall consider that
\begin{equation}\label{Ups}
\mathbf{\Upsilon}_{\epsilon_m}:=\mathbf{\Upsilon}^1_{\epsilon_m}+\mathbf{\Upsilon}^2_{\epsilon_m}
\end{equation}
is extended from $\omega$ to $\mathds{R}^{N}$ by zero, where
\begin{equation}\label{Ups-1}
\mathbf{\Upsilon}^1_{\epsilon_m}:=-\left(|\mathbf{D}(\mathbf{u}_{\epsilon_m})|^{\gamma-2}\mathbf{D}(\mathbf{u}_{\epsilon_m})-\mathbf{S}_1+|\mathbf{D}(\mathbf{u}_{\epsilon_m})|^{q(\mathbf{x})-2}\mathbf{D}(\mathbf{u}_{\epsilon_m})-\mathbf{S}_2\right)+p^1_{\epsilon_m}\mathbf{I},
\end{equation}
\begin{equation}\label{Ups-2}
\mathbf{\Upsilon}^2_{\epsilon_m}:=\mathbf{u}_{\epsilon_m}\otimes\mathbf{u}_{\epsilon_m}\Phi_{\epsilon_m}(|\mathbf{u}_{\epsilon_m}|)-                                    \mathbf{u}\otimes\mathbf{u}+p^2_{\epsilon_m}\mathbf{I},
\end{equation}
and $\mathbf{I}$ denotes the identity tensor.

Now, due to the definition (\ref{v-em-chi}) and by virtue of (\ref{convg-La}) and (\ref{convg-s-g}), we have
\begin{equation}\label{un-b-v-em}
\mbox{$\mathbf{w}_{\epsilon_m}\to \mathbf{0}$\quad weakly in $\mathbf{W}^{1,{\gamma}}(\mathds{R}^{N})$,\quad as $m\to\infty$,}
\end{equation}
\begin{equation}\label{str-v-em}
\mbox{$\mathbf{w}_{\epsilon_m}\to \mathbf{0}$\quad strongly in $\mathbf{L}^{\kappa}(\mathds{R}^{N})$,\quad as $m\to\infty$,
\quad
for any $\kappa:1\leq\kappa<\gamma^{\ast}$.
}
\end{equation}
Moreover, due to (\ref{conv-Re-0}), (\ref{convg-S-b'}) and (\ref{bd2-p1-em}) by one hand, and due to (\ref{convg-Phiuxu}), (\ref{limit-eq-G}) and (\ref{bd2-p2-em}) on the other, we have
\begin{equation}\label{un-b-T-p1-em}
\|\mathbf{\Upsilon}^1_{\epsilon_m}\|_{\mathbf{L}^{\gamma'}(\mathds{R}^{N})}\leq C,
\end{equation}
\begin{equation}\label{un-bA-uxu-p2-em}
\|\mathbf{\Upsilon}^2_{\epsilon_m}\|_{\mathbf{L}^{\frac{\gamma^{\ast}}{2}}(\mathds{R}^{N})}\leq C.
\end{equation}
In addition to (\ref{un-bA-uxu-p2-em}), we see that, due to (\ref{convg-s-g}) and (\ref{bd2-p2-em}),
\begin{equation}\label{un-b-uxu-p2-em}
\mathbf{\Upsilon}^2_{\epsilon_m}\to 0\quad\mbox{strongly in $\mathbf{L}^{\frac{\kappa}{2}}(\mathds{R}^{N})$,\quad as $m\to\infty$,
\quad
for any $\kappa:1\leq\kappa<\gamma^{\ast}$.}
\end{equation}

Next, let us consider the Hardy-Littlewood maximal functions of $|\mathbf{w}_{\epsilon_m}|$ and $|\mathbf{\nabla}\mathbf{w}_{\epsilon_m}|$ defined by
\begin{eqnarray*}
  &\displaystyle& \mathcal{M}(|\mathbf{w}_{\epsilon_m}|)(\mathbf{x}):=\sup_{0<R<\infty}\frac{1}{\mathcal{L}_N(B_R(\mathbf{x}))}\int_{B_R(\mathbf{x})}|\mathbf{w}_{\epsilon_m}(\mathbf{y})|\,d\mathbf{y},\\
   &\displaystyle& \mathcal{M}(|\mathbf{\nabla}\mathbf{w}_{\epsilon_m}|)(\mathbf{x}):=\sup_{0<R<\infty}\frac{1}{\mathcal{L}_N(B_R(\mathbf{x}))}\int_{B_R(\mathbf{x})}|\mathbf{\nabla}\mathbf{w}_{\epsilon_m}(\mathbf{y})|\,d\mathbf{y};
\end{eqnarray*}
where $B_R(\mathbf{x})$ denotes the ball of $\mathds{R}^N$ centered at $\mathbf{x}$ and with radius $R>0$, and $\mathcal{L}_N(\omega)$ is the $N$-dimensional Lebesgue measure of $\omega$. %
Arguing as in Diening \emph{et al.}~\cite[p. 218]{DMS-2008} and using the boundedness of the Hardy-Littlewood maximal operator $\mathcal{M}$ (see \emph{e.g.} Stein~\cite[Theorem I.1.1]{Stein-1970}), we can prove that for all $m\in\mathds{N}$ and all $j\in\mathds{N}_0$ there exists
\begin{equation}\label{int-2j}
\lambda_{m,j}\in\left[2^{2^j},2^{2^{j+1}}\right)
\end{equation}
such that
\begin{equation}\label{Leb-F}
\mathcal{L}_{N}\left(F_{m,j}\right)\leq
2^{-j}\lambda_{m,j}^{-\kappa}\,\|\mathbf{w}_{\epsilon_m}\|_{\mathbf{L}^{\kappa}(\mathds{R}^{N})},
\quad
\mbox{for any $\kappa:1\leq\kappa<\gamma^{\ast}$,}
\end{equation}
\begin{equation}\label{Leb-G}
\mathcal{L}_{N}\left(G_{m,j}\right)\leq
2^{-j}\lambda_{m,j}^{-\gamma}\,\|\mathbf{\nabla}\mathbf{w}_{\epsilon_m}\|_{\mathbf{L}^{\gamma}(\mathds{R}^{N})},
\end{equation}
where
\begin{eqnarray*}
  &\displaystyle& F_{m,j}:=\left\{\mathbf{x}\in\mathds{R}^{N}:\mathcal{M}(|\mathbf{w}_{\epsilon_m}|)(\mathbf{x})>2\lambda_{m,j}\right\},\\
   &\displaystyle& G_{m,j}:=\left\{\mathbf{x}\in\mathds{R}^{N}:\mathcal{M}(|\mathbf{\nabla}\mathbf{w}_{\epsilon_m}|)(\mathbf{x})>2\lambda_{m,j}\right\}.
\end{eqnarray*}
Setting
\begin{equation}\label{R-mj}
R_{m,j}:=F_{m,j}\cup G_{m,j}\cup\left\{\mathbf{x}\in\mathds{R}^{N}:\ \mbox{$\mathbf{x}$ is not a Lebesgue point of $|\mathbf{w}_{\epsilon_m}|$}\right\},
\end{equation}
we can see that, by virtue of (\ref{Leb-F})-(\ref{R-mj}) and (\ref{un-b-v-em})-(\ref{str-v-em}),
\begin{equation}\label{lim-m-R}
\limsup_{m\to\infty}\mathcal{L}_{N}\left(R_{m,j}\right)\leq \limsup_{m\to\infty}C2^{-j}\lambda_{m,j}^{-\gamma}.
\end{equation}

The following result helps us to approximate $W^{1,\gamma}$-functions by Lipschitz ones.

\begin{lemma}\label{AF-Lem}
Let $\omega\subset\mathds{R}^N$ be an open bounded set with a Lipschitz-continuous boundary $\partial\omega$.
Assume that $s\geq 1$ and let $\mathbf{u}\in\mathbf{W}^{1,s}(\omega)$. %
Then, for every $\lambda_1,\ \lambda_2>0$, there exists $\mathbf{v}\in\mathbf{W}^{1,\infty}(\omega)$ such that
\begin{equation*}
\|\mathbf{v}\|_{\mathbf{L}^{\infty}(\omega)}\leq\lambda_1,
\end{equation*}
\begin{equation*}
\|\mathbf{\nabla}\mathbf{v}\|_{\mathbf{L}^{\infty}(\omega)}\leq C\lambda_2,\quad C=C(N,\omega).
\end{equation*}
In addition,
\begin{equation*}
\begin{split}
  & \{\mathbf{x}\in\omega:\mathbf{v}(\mathbf{x})\not=\mathbf{u}(\mathbf{x})\}\subset \\
  &\ \omega\cap
\left(\{\mathbf{x}\in\mathds{R}^N:\mathcal{M}(|\mathbf{u}|)(\mathbf{x})>\lambda_1\}\cup
\{\mathbf{x}\in\mathds{R}^N:\mathcal{M}(|\mathbf{\nabla}\mathbf{u}|)(\mathbf{x})>\lambda_2\}\cup\right.\\
  &\ \qquad\left.
\{\mathbf{x}\in\mathds{R}^{N}:\ \mbox{$\mathbf{x}$ is not a Lebesgue point of $|\mathbf{u}|$}\}\right).
\end{split}
\end{equation*}
\end{lemma}
\begin{proof}
See Acerbi and Fusco~\cite{AF-1988} for the original result and Landes~\cite[Proposition 2.2]{Landes-1996}
for the last statement.
\end{proof}

\noindent Then, by Lemma~\ref{AF-Lem} together with the definition of $\mathbf{w}_{\epsilon_m}$ (\emph{f.} (\ref{v-em-chi})), there exists
\begin{equation}\label{z-mj}
\mathbf{z}_{m,j}\in\mathbf{W}^{1,\infty}(\mathds{R}^N),\qquad
\mathbf{z}_{m,j}=\left\{\begin{array}{cc}
                   \mathbf{w}_{\epsilon_m} & \mbox{in $\omega\setminus A_{m,j}$} \\
                   0 & \mathds{R}^N\setminus\omega
                 \end{array}\right.\,,
\end{equation}
where
\begin{equation}\label{A-mj}
A_{m,j}:=\{\mathbf{x}\in\omega:\mathbf{z}_{m,j}(\mathbf{x})\not=\mathbf{w}_{\epsilon_m}(\mathbf{x})\},
\end{equation}
such that
\begin{equation}\label{eq1-AF}
\|\mathbf{z}_{m,j}\|_{\mathbf{L}^{\infty}(\omega)}\leq2\lambda_{m,j},
\end{equation}
\begin{equation}\label{eq2-AF}
\|\mathbf{\nabla}\mathbf{z}_{m,j}\|_{\mathbf{L}^{\infty}(\omega)}\leq C\lambda_{m,j},\quad C=C(N,\omega).
\end{equation}
Moreover, by the last statement of the above lemma and using the notations (\ref{Leb-F})-(\ref{R-mj}) and (\ref{A-mj}),
\begin{equation}\label{eq3-AF}
A_{m,j}\subset\omega\cap R_{m,j}.
\end{equation}
As a consequence of (\ref{lim-m-R}) and (\ref{eq3-AF}),
\begin{equation}\label{lim-m-A}
\limsup_{m\to\infty}\mathcal{L}_{N}\left(A_{m,j}\right)\leq \limsup_{m\to\infty}
C2^{-j}\lambda_{m,j}^{-\gamma}.
\end{equation}
On the other hand, due to (\ref{un-b-v-em}), (\ref{eq1-AF})-(\ref{eq2-AF}) and (\ref{lim-m-A}), we can prove that for any $j\in\mathds{N}_0$
\begin{equation}\label{w-z-em-a}
\mbox{$\mathbf{z}_{m,j}\to \mathbf{0}$\quad weakly in $\mathbf{W}^{1,{\gamma}}_0(\omega)$,\quad as $m\to\infty$.}
\end{equation}
Then by Sobolev's compact imbedding theorem, we get for any $j\in\mathds{N}_0$
\begin{equation*}
\mbox{$\mathbf{z}_{m,j}\to \mathbf{0}$ strongly in $\mathbf{L}^{\kappa}(\omega)$,\quad as $m\to\infty$,\quad for any $\kappa: 1\leq\kappa<\gamma^{\ast}$.}
\end{equation*}
Using this information, (\ref{eq1-AF}) and interpolation, we prove that for any $j\in\mathds{N}_0$
\begin{equation}\label{str-z-em}
\mbox{$\mathbf{z}_{m,j}\to \mathbf{0}$\quad strongly in $\mathbf{L}^s(\omega)$,\quad as $m\to\infty$,\quad for any $s:1\leq s<\infty$.}
\end{equation}
Finally, as a consequence of (\ref{w-z-em-a}) and (\ref{str-z-em}), we obtain for any $j\in\mathds{N}_0$
\begin{equation}\label{w-con-s-z-em}
\mbox{$\mathbf{z}_{m,j}\to \mathbf{0}$\quad weakly in $\mathbf{W}^{1,s}_0(\omega)$,\quad as $m\to\infty$,\quad for any $s:1\leq s<\infty$.}
\end{equation}

\section{Convergence of the approximated extra stress tensor}\label{Sect-Conv-EST}

Let us first observe that, using the notations (\ref{Ups})-(\ref{Ups-2}), we can write (\ref{dist2-der}) in the following simplified form
\begin{equation}\label{dist3-der}
\mathbf{div}\mathbf{\Upsilon}_{\epsilon_m}=\mathbf{0}
\quad\mbox{in\quad $\mathcal{D}'(\omega)$}.
\end{equation}
On the other hand, due to (\ref{un-b-T-p1-em})-(\ref{un-bA-uxu-p2-em}), $\mathbf{\Upsilon}_{\epsilon_m}\in\mathbf{L}^{r}(\mathds{R}^N)$ for $r$ satisfying to (\ref{r}). %
Then, using this information and (\ref{w-con-s-z-em}), we infer, from (\ref{dist3-der}), that for any $j\in\mathds{N}_0$
\begin{equation}\label{dist4-der}
  \int_{\omega}\mathbf{\Upsilon}_{\epsilon_m}:\mathbf{\nabla}\mathbf{z}_{m,j}\,d\mathbf{x}=0.
\end{equation}
Expanding $\mathbf{\Upsilon}_{\epsilon_m}$ in (\ref{dist4-der}) through the notations (\ref{Ups})-(\ref{Ups-2}) and subtracting and adding the integral
\begin{equation*}
\int_{\omega}\left(\mu_1|\mathbf{D}(\mathbf{u})|^{\gamma-2}+\mu_2|\mathbf{D}(\mathbf{u})|^{q(\mathbf{x})-2}\right)
\mathbf{D}(\mathbf{u}):\mathbf{D}(\mathbf{z}_{m,j})\,d\mathbf{x}
\end{equation*}
to the left hand side of the resulting equation, we obtain for any $j\in\mathds{N}_0$
\begin{equation}\label{dist5-der}
\begin{split}
  & \mu_1\int_{\omega}\left(|\mathbf{D}(\mathbf{u}_{\epsilon_m})|^{\gamma-2}\mathbf{D}(\mathbf{u}_{\epsilon_m})-|\mathbf{D}(\mathbf{u})|^{\gamma-2}\mathbf{D}(\mathbf{u})\right):\mathbf{D}(\mathbf{z}_{m,j})\,d\mathbf{x} +\\
  &\mu_2\int_{\omega}\left(|\mathbf{D}(\mathbf{u}_{\epsilon_m})|^{q(\mathbf{x})-2}\mathbf{D}(\mathbf{u}_{\epsilon_m})-|\mathbf{D}(\mathbf{u})|^{q(\mathbf{x})-2}\mathbf{D}(\mathbf{u})\right):\mathbf{D}(\mathbf{z}_{m,j})\,d\mathbf{x}  =\\
  &\mu_1\int_{\omega}\left(\mathbf{S}_1-|\mathbf{D}(\mathbf{u})|^{\gamma-2}\mathbf{D}(\mathbf{u})\right):\mathbf{D}(\mathbf{z}_{m,j})\,d\mathbf{x}+ \\
  &\mu_2\int_{\omega}\left(\mathbf{S}_2-|\mathbf{D}(\mathbf{u})|^{q(\mathbf{x})-2}\mathbf{D}(\mathbf{u})\right):\mathbf{D}(\mathbf{z}_{m,j})\,d\mathbf{x}  +\\
  &\int_{\omega}p^1_{\epsilon_m}\,\mathrm{div}\,\mathbf{z}_{m,j}\,d\mathbf{x}+\\
  &\int_{\omega}\left(\mathbf{u}_{\epsilon_m}\otimes\mathbf{u}_{\epsilon_m}\Phi_{\epsilon_m}(|\mathbf{u}_{\epsilon_m}|)-                                    \mathbf{u}\otimes\mathbf{u}+p^2_{\epsilon_m}\mathbf{I}\right):\mathbf{D}(\mathbf{z}_{m,j})\,d\mathbf{x}\\
  &:=J_{m,j}^1+J_{m,j}^2+J_{m,j}^3+J_{m,j}^4.
  \end{split}
\end{equation}
We claim that, for a fixed $j$,
\begin{equation}\label{dist6-der}
\begin{split}
  &\mu_1\lim_{m\to\infty}\int_{\omega}\left(|\mathbf{D}(\mathbf{u}_{\epsilon_m})|^{\gamma-2}\mathbf{D}(\mathbf{u}_{\epsilon_m})-|\mathbf{D}(\mathbf{u})|^{\gamma-2}\mathbf{D}(\mathbf{u})\right):\mathbf{D}(\mathbf{z}_{m,j})\,d\mathbf{x} +\\
  &\mu_2\lim_{m\to\infty}\int_{\omega}\left(|\mathbf{D}(\mathbf{u}_{\epsilon_m})|^{q(\mathbf{x})-2}\mathbf{D}(\mathbf{u}_{\epsilon_m})-|\mathbf{D}(\mathbf{u})|^{q(\mathbf{x})-2}\mathbf{D}(\mathbf{u})\right):\mathbf{D}(\mathbf{z}_{m,j})\,d\mathbf{x} \leq C2^{-\frac{j}{\gamma}}
\end{split}
\end{equation}
To prove this, we will carry out the passage to the limit $m\to\infty$ in all absolute values $|J_{m,j}^i|$, $i=1,\dots,4$.

$\bullet\ \limsup_{m\to\infty}|J_{m,j}^1|=0$. By (\ref{w-con-s-z-em}), with $s=\gamma$, this is true once we can justify that
$\mathbf{S}_1-|\mathbf{D}(\mathbf{u})|^{\gamma-2}\mathbf{D}(\mathbf{u})$ is uniformly bounded in $\mathbf{L}^{\gamma'}(\omega)$. %
But this is an immediate consequence of (\ref{conv-Re-0}).

$\bullet\ \limsup_{m\to\infty}|J_{m,j}^2|=0$. Again, by (\ref{w-con-s-z-em}) with $s=\gamma$, this is true if
$\mathbf{S}_2-|\mathbf{D}(\mathbf{u})|^{q(\mathbf{x})-2}\mathbf{D}(\mathbf{u})$ is uniformly bounded in $\mathbf{L}^{\gamma'}(\omega)$. %
Analogously to the previous case, we can justify this now by using (\ref{convg-S-b'}).

$\bullet\ \limsup_{m\to\infty}|J_{m,j}^3|\leq C 2^{-\frac{j}{\gamma}}$. In fact, by Hölder's inequality and (\ref{bd2-p1-em}) together with (\ref{conv-Re-0}) and (\ref{convg-S-b'}), and using the definition of $\mathbf{z}_{m,j}$ (\emph{f.} (\ref{z-mj})) together with the fact that $\mathrm{div}\mathbf{w}_{\epsilon_m}=0$ in $\omega$ (\emph{f.} (\ref{v-em-chi})),
\begin{equation*}
|J_{m,j}^3|\leq C_1\|\mathrm{div}\,\mathbf{z}_{m,j}\|_{\mathbf{L}^{\beta}(\omega)}
     \leq C_1\|\mathbf{\nabla}\mathbf{z}_{m,j}\|_{\mathbf{L}^{\infty}(\omega)}\mathcal{L}_N(A_{m,j})^{\frac{1}{\gamma}}.
\end{equation*}
Then, by the application of (\ref{eq2-AF}) and (\ref{lim-m-A}), it follows
\begin{equation}\label{cond1-HI}
|J_{m,j}^3|\leq C_2\limsup_{m\to\infty}\lambda_{m,j}\left(2^{-j}\lambda^{-\gamma}\right)^{\frac{1}{\gamma}}\leq
C_2 2^{-\frac{j}{\gamma}}.
\end{equation}

$\bullet\ \limsup_{m\to\infty}|J_{m,j}^4|=0$. Using Hölder's inequality and the notation (\ref{Ups-2}), we have
\begin{equation*}
\begin{split}
  |J_{m,j}^4|&\leq \|\mathbf{\Upsilon}^2_{\epsilon_m}\|_{\mathbf{L}^{1}(\omega)}\|\mathbf{\nabla}\mathbf{z}_{m,j}\|_{\mathbf{L}^{\infty}(\omega)} \\
       &\leq C_1\|\mathbf{\Upsilon}^2_{\epsilon_m}\|_{\mathbf{L}^{1}(\omega)}\to 0,\quad\mbox{as}\ m\to\infty.
\end{split}
\end{equation*}
The last inequality and the conclusion follow, respectively, from (\ref{eq2-AF}) and  (\ref{un-b-uxu-p2-em}) with $\kappa=2$, observing that here the assumption (\ref{al-TE-1}) implies $2<\gamma^{\ast}$.

Gathering the estimates above we just have proven (\ref{dist6-der}).

We proceed with the proof by using an argument due to Dal Maso and Murat~\cite[Theorem 5]{DalM-M} (see also the references cited in \cite[Remark 4]{DalM-M}).
Firstly, observing the definition of $\mathbf{z}_{m,j}$ (\emph{f.} (\ref{z-mj})), we have
\begin{equation}\label{eq21-est-theta}
\begin{split}
&\mu_1\int_{\omega}\left(|\mathbf{D}(\mathbf{u}_{\epsilon_m})|^{\gamma-2}\mathbf{D}(\mathbf{u}_{\epsilon_m})-|\mathbf{D}(\mathbf{u})|^{\gamma-2}\mathbf{D}(\mathbf{u})\right):\mathbf{D}(\mathbf{z}_{m,j})\,d\mathbf{x} +\\
&\mu_2\int_{\omega}\left(|\mathbf{D}(\mathbf{u}_{\epsilon_m})|^{q(\mathbf{x})-2}\mathbf{D}(\mathbf{u}_{\epsilon_m})-|\mathbf{D}(\mathbf{u})|^{q(\mathbf{x})-2}\mathbf{D}(\mathbf{u})\right):\mathbf{D}(\mathbf{z}_{m,j})\,d\mathbf{x}  \\&:=
\mu_1I_{m,j}^1+\mu_2I_{m,j}^2+\mu_1II_{m,j}^1+\mu_2II_{m,j}^2,
\end{split}
\end{equation}
where
\begin{equation*}
\begin{split}
&I_{m,j}^1:=\int_{\omega\setminus A_{m,j}}\left(|\mathbf{D}(\mathbf{u}_{\epsilon_m})|^{\gamma-2}\mathbf{D}(\mathbf{u}_{\epsilon_m})-|\mathbf{D}(\mathbf{u})|^{\gamma-2}\mathbf{D}(\mathbf{u})\right):(\mathbf{D}(\mathbf{u}_{\epsilon_m})-\mathbf{D}(\mathbf{u}))\,d\mathbf{x},\\
&I_{m,j}^2:=\int_{\omega\setminus A_{m,j}}\left(|\mathbf{D}(\mathbf{u}_{\epsilon_m})|^{q(\mathbf{x})-2}\mathbf{D}(\mathbf{u}_{\epsilon_m})-|\mathbf{D}(\mathbf{u})|^{q(\mathbf{x})-2}\mathbf{D}(\mathbf{u})\right):(\mathbf{D}(\mathbf{u}_{\epsilon_m})-\mathbf{D}(\mathbf{u}))\,d\mathbf{x},\\
&II_{m,j}^1:= \int_{A_{m,j}}\left(|\mathbf{D}(\mathbf{u}_{\epsilon_m})|^{\gamma-2}\mathbf{D}(\mathbf{u}_{\epsilon_m})-|\mathbf{D}(\mathbf{u})|^{\gamma-2}\mathbf{D}(\mathbf{u})\right):\mathbf{D}(\mathbf{z}_{m,j})\,d\mathbf{x},\\
&II_{m,j}^2:= \int_{A_{m,j}}\left(|\mathbf{D}(\mathbf{u}_{\epsilon_m})|^{q(\mathbf{x})-2}\mathbf{D}(\mathbf{u}_{\epsilon_m})-|\mathbf{D}(\mathbf{u})|^{q(\mathbf{x})-2}\mathbf{D}(\mathbf{u})\right):\mathbf{D}(\mathbf{z}_{m,j})\,d\mathbf{x}.
     \end{split}
\end{equation*}
Then (\ref{dist6-der}) and (\ref{eq21-est-theta}) imply that
\begin{equation}\label{est-I1}
\limsup_{m\to\infty}\left(I_{m,j}^1+I_{m,j}^2\right)\leq C\left[\limsup_{m\to\infty}\left(|II_{m,j}^1+II_{m,j}^2|\right)+2^{-\frac{j}{\gamma}}\right].
\end{equation}
For the terms $II_{m,j}^1$ and $II_{m,j}^2$, we have by applying successively Hölder's inequality, (\ref{conv-Re-0}) and (\ref{convg-S-b'}), altogether with (\ref{eq2-AF}),
\begin{equation*}
\begin{split}
&|II_{m,j}^1+II_{m,j}^2|\leq\\
&C_1
\|\left(|\mathbf{D}(\mathbf{u}_{\epsilon_m})|^{\gamma-2}\mathbf{D}(\mathbf{u}_{\epsilon_m})-|\mathbf{D}(\mathbf{u})|^{\gamma-2}\mathbf{D}(\mathbf{u})\right)
  \|_{\mathbf{L}^{\gamma'}(A_{m,j})}\|\mathbf{\nabla}\mathbf{z}_{m,j}\|_{\mathbf{L}^{\gamma}(A_{m,j})}+\\
&C_2
\|\left(|\mathbf{D}(\mathbf{u}_{\epsilon_m})|^{q(\mathbf{x})-2}\mathbf{D}(\mathbf{u}_{\epsilon_m})-|\mathbf{D}(\mathbf{u})|^{q(\mathbf{x})-2}\mathbf{D}(\mathbf{u})\right)
  \|_{\mathbf{L}^{\gamma'}(A_{m,j})}
\|\mathbf{\nabla}\mathbf{z}_{m,j}\|_{\mathbf{L}^{\gamma}(A_{m,j})}\leq\\
&
C_3\lambda_{m,j}\mathcal{L}_N(A_{m,j})^{\frac{1}{\gamma}}.
\end{split}
\end{equation*}
Then using the same arguing as we did for (\ref{cond1-HI}) it holds for any $j\in\mathds{N}_0$
\begin{equation}\label{est-I2}
\limsup_{m\to\infty}|II_{m,j}^1+II_{m,j}^2|\leq
C2^{-\frac{j}{\gamma}}.
\end{equation}
As a consequence of (\ref{est-I1}) and (\ref{est-I2}), we obtain for any $j\in\mathds{N}_0$
\begin{equation}\label{est2-I1}
\limsup_{m\to\infty}\left(I_{m,j}^1+I_{m,j}^2\right)\leq
C2^{-\frac{j}{\gamma}}.
\end{equation}
Now, by Hölder's inequality and having in mind the definition of $\mathbf{z}_{m,j}$ (\emph{f.} (\ref{z-mj})), we have for any $\theta\in(0,1)$
\begin{equation}\label{est-g-theta}
\begin{split}
\int_{\omega}g_{\epsilon_m}^{\theta}\,d\mathbf{x}
\leq &C_1\left[(\mathcal{I}_{m,j}^1)^{\theta}\mathcal{L}_N(\omega\setminus A_{m,j})^{1-\theta}+
(\mathcal{II}_{m,j}^1)^{\theta}\mathcal{L}_N(A_{m,j})^{1-\theta}\right]+\\
&C_2\left[(\mathcal{I}_{m,j}^2)^{\theta}\mathcal{L}_N(\omega\setminus A_{m,j})^{1-\theta}+
(\mathcal{II}_{m,j}^2)^{\theta}\mathcal{L}_N(A_{m,j})^{1-\theta}\right]
\end{split}
\end{equation}
where
\begin{equation*}
\begin{split}
  g_{\epsilon_m}:=&\mu_1\left||\mathbf{D}(\mathbf{u}_{\epsilon_m})|^{\gamma-2}\mathbf{D}(\mathbf{u}_{\epsilon_m})-|\mathbf{D}(\mathbf{u})|^{\gamma-2}\mathbf{D}(\mathbf{u})\right|+\\
                  &\mu_2\left||\mathbf{D}(\mathbf{u}_{\epsilon_m})|^{q(\mathbf{x})-2}\mathbf{D}(\mathbf{u}_{\epsilon_m})-|\mathbf{D}(\mathbf{u})|^{q(\mathbf{x})-2}\mathbf{D}(\mathbf{u})\right|,
\end{split}
\end{equation*}
and $\mathcal{I}_{m,j}^1$, $\mathcal{II}_{m,j}^1$, $\mathcal{I}_{m,j}^2$ and $\mathcal{II}_{m,j}^2$ are, respectively, the functions $I_{m,j}^1$, $II_{m,j}^1$, $I_{m,j}^2$ and $II_{m,j}^2$, but with the integrand functions replaced by their absolute values.
Arguing as we did to prove (\ref{est-I2})-(\ref{est2-I1}) and using (\ref{lim-m-A}), it follows from (\ref{est-g-theta}) that
\begin{equation}\label{eq3-est-theta}
\limsup_{m\to\infty}\int_{\omega}g_{\epsilon_m}^{\theta}\,d\mathbf{x}
\leq
C_12^{-\theta\frac{j}{\gamma}}+C_22^{-\theta\frac{j}{\gamma}-(1-\theta)j}\,.
\end{equation}
Since $\gamma>1$, $\theta\in(0,1)$ and $j\in\mathds{N}_0$ is arbitrary, $2^{-\theta\frac{j}{\gamma}}\to 0$ and $2^{-\theta\frac{j}{\gamma}-(1-\theta)j}\to 0$, as $j\to\infty$.
This and (\ref{eq3-est-theta}) imply that for any $\theta\in(0,1)$
\begin{equation*}
\limsup_{m\to\infty}\int_{\omega}g_{\epsilon_m}^{\theta}\,d\mathbf{x}=0.
\end{equation*}
Then, passing to a subsequence,
\begin{equation}\label{eq2-est-theta}
g_{\epsilon_m}\to 0\quad\mbox{a.e. in $\omega$},\quad\mbox{as $m\to\infty$}.
\end{equation}
Due to the fact that $|\mathbf{D}(\mathbf{u})|^{\gamma-2}\mathbf{D}(\mathbf{u})$ and $|\mathbf{D}(\mathbf{u})|^{q(\mathbf{x})-2}\mathbf{D}(\mathbf{u})$ are strictly monotonous and continuous tensors on $\mathbf{D}(\mathbf{u})$, we can apply  \cite[Lemma 6]{DalM-M} together with (\ref{eq2-est-theta}) (see also Lions~\cite[Lemme 2.2.2]{Lions-1969}), to establish that
\begin{equation}\label{eq4-est-theta}
\mathbf{D}(\mathbf{u}_{\epsilon_m})\to\mathbf{D}(\mathbf{u})\quad\mbox{a.e. in $\omega$},\quad\mbox{as $m\to\infty$}.
\end{equation}
Finally, (\ref{est-gam'})-(\ref{est-q'}) and (\ref{eq4-est-theta}) allow us to use Vitali's theorem together with (\ref{conv-Re-0})-(\ref{convg-S-b'}) to conclude that $\mathbf{S}_1=|\mathbf{D}(\mathbf{u})|^{\gamma-2}\mathbf{D}(\mathbf{u})$ and $\mathbf{S}_2=|\mathbf{D}(\mathbf{u})|^{q(\mathbf{x})-2}\mathbf{D}(\mathbf{u})$.

This concludes the proof of Theorem~\ref{th-exst-ws-pp}. $\blacksquare$

\begin{remark}
Using the parabolic versions of the techniques we have used here (see \emph{e.g.} Oliveira~\cite{O-2012}), we can establish an existence result for the parabolic model of the problem (\ref{geq1-inc})-(\ref{e-stress}). This work will be the aim of a forthcoming paper.
\end{remark}

\end{document}